\documentclass{article}
\usepackage[utf8]{inputenc}
\usepackage{lineno,epsfig,graphics,amssymb,amsmath,graphicx,pifont,soul,esint,float,multirow,array,color,url,booktabs,mathrsfs,mathdots,amsthm}
\usepackage[dvipsnames]{xcolor}
\usepackage{siunitx}
\usepackage{enumitem}
\usepackage[margin=1in]{geometry}
\usepackage{url}
\usepackage{caption}
\usepackage{subcaption}
\usepackage{cleveref}
\usepackage{mathtools}
\usepackage{authblk}

\newcommand{\bl}[1]{\boldsymbol{#1}}
\newcommand{\dnorm}[1]{{\left\vert\kern-0.25ex\left\vert\kern-0.25ex\left\vert #1 
    \right\vert\kern-0.25ex\right\vert\kern-0.25ex\right\vert}}

\newcommand{\inner}[1]{\bigg({#1}\bigg)}
\newcommand{\Gh}[0]{{\Gamma_{\rm h}}}
\newcommand{\Gg}[0]{{\Gamma_{\rm g}}}
\newcommand{\Ot}[0]{{\tilde \Omega}}

\newcommand{\res}[0]{{\pmb{\mathscr{L}}}}

\title{Introducing a Harmonic Balance Navier-Stokes Finite Element Solver to Accelerate Cardiovascular Simulations}
\author[1]{Dongjie Jia\thanks{Corresponding Author: jia191@purdue.edu. Present address: Purdue University, West Lafayette, IN, USA.}}
\author[1]{Mahdi Esmaily\thanks: }
\affil[1]{Sibley School of Mechanical and Aerospace Engineering, Cornell University, Ithaca, NY, USA}

\begin{document}

\maketitle

\begin{abstract}
    The adoption of cardiovascular simulations for diagnosis and surgical planning on a patient-specific basis requires the development of faster methods than the existing state-of-the-art techniques. To address this need, we leverage the periodic nature of these flows to accurately capture their time-dependence using spectral discretization. Owing to the reduced size of the discrete problem, the resulting approach, known as the harmonic balance method, significantly lowers the solution cost when compared against the conventional time marching methods. This study describes a stabilized finite element implementation of the harmonic balanced method that targets the simulation of physically-stable time-periodic flows. That stabilized method is based on the Galerkin/least-squares formulation that permits stable solution in convection-dominant flows and convenient use of the same interpolation functions for velocity and pressure. We test this solver against its equivalent time marching method using three common physiological cases where blood flow is modeled in a Glenn operation, a cerebral artery, and a left main coronary artery. Using the conventional time marching solver, simulating these cases takes more than ten hours. That cost is reduced by up to two orders of magnitude when the proposed harmonic balance solver is utilized, where a solution is produced in approximately 30 minutes. We show that that solution is in excellent agreement with the conventional solvers when the number of modes is sufficiently large (about 10) to accurately represent the imposed boundary conditions.
\end{abstract}

\section{Introduction}

Computational Fluid Dynamics (CFD) has played a foundational role in studying cardiovascular blood flows in the recent decades~\cite{taylor1998finite,figueroa2006coupled,moghadam2013modular,morris2016computational}. 
These simulations provide a non-invasive way of accessing information that is difficult or impossible to obtain using conventional clinical methods. 
Simulation results such as flow paths, pressure distribution, and wall shear stress can provide critical information for non-invasive diagnosis~\cite{antiga2008image,kim2010patient,chung2015cfd,mittal2016computational}, predictive surgical planning~\cite{pekkan2008patient,esmaily2015simulations,jia2019image}, and surgical optimization~\cite{esmaily2012optimization, marsden2014optimization,jia2022characterization}. 
Unfortunately, the applications are mainly limited to research settings due to the cost of performing high-accuracy simulations, which ranges from hours to days when run on a dedicated computing cluster~\cite{trusty2018fontan,jia2021efficient}. 
Such cost can become a limiting factor in real-life patient-specific settings considering the availability of computational resources and the short time frame available for clinical decision-making.

One current cost-saving strategy for CFD simulation in clinical settings involves assuming steady-state flows. 
Simulating transient cardiovascular flows is a resource-intensive process, which requires repeatedly solving millions of degrees of freedom for thousands of time steps each cardiac cycle for multiple cycles to reach solution convergence. 
When the steady-state assumption is used instead of transient analysis, the required number of time steps is significantly reduced, which in turn reduces the simulation cost. 
An example in practice is calculating fractional flow reserve (FFR) to analyze coronary artery disease, which takes minutes to compute the steady-state CFD result~\cite{taylor2013computational,morris2015virtual,driessen2019comparison,modi2019predicting}. 
However, because the steady-state assumption does not capture the cardiovascular system's dynamic and pulsatile behavior, the results carry significant errors compared to measured FFR values. 
Furthermore, steady-state CFD cannot provide estimations for transient-dependent parameters, such as maximum wall shear stress and pressure.
Other strategies use modeling to circumvent performing CFD in clinical settings and eliminate the computation cost associated with simulations.
These modeling techniques, such as reduced order modeling~\cite{liang2005closed,mynard2012simple,mirramezani2019reduced,pfaller2022automated}, and more recently, data-driven machine learning~\cite{tesche2018coronary,li2021prediction,schwarz2023beyond}, uses prior obtained CFD results or clinical data to construct data-informed models.
Although these models are cost-efficient for clinical use, they struggle to accommodate patient-specific variations in conditions and geometries. 
In addition, these models do not yield high resolution three-dimensional temporal-spacial results like CFD.
A majority of these models also do not fundamentally satisfy fluid conservation laws, making their hemodynamic predictions error-prone. 
Therefore, transient CFD remains the gold standard for accurate patient-specific cardiovascular simulations.

As mentioned above, conventional transient CFD simulation uses time-stepping to advance the solution in time. 
The time-stepping scheme is fundamentally responsible for the high simulation cost due to the significant number of time steps required per cardiac cycle and multiple cycles to achieve converged periodic results. 
On the other hand, the periodic and smooth temporal behavior of physically stable cardiovascular flows can be well-approximated using a few frequency modes. 
Solving those selective frequency modes can reduce the simulation cost significantly compared to the conventional time-stepping scheme.
In a recent series of studies, we suggested conducting transient cardiovascular CFD simulations in the time-spectral (or frequency) domain~\cite{meng2021scalable,esmaily2023stabilized,esmaily2025augmented,esmaily2024new}.
The implementation of the frequency finite element solver requires several innovative approaches, starting with a specialized code structure to separate the complex values to use only real arithmetics in the solver and a pressure stabilization scheme to allow equal shape functions for pressure and velocity~\cite{meng2021scalable,esmaily2023stabilized}. 
In addition, a new finite element stabilization method was proposed to improve solution accuracy and compared against other existing stabilization methods in solving the frequency convection-diffusion equation~\cite{esmaily2025augmented}.

The stabilized finite element method is widely used for cardiovascular simulations due owing to its ability to handle complex geometries, unstructured grid, and fluid-structure interaction~\cite{van2003finite,bazilevs2009computational,quarteroni2016geometric,kamensky2017immersogeometric}.
In our latest published work, we developed a CFD solver that uses a specialized stabilized finite element method to solve the Navier-Stokes equations in frequency domain~\cite{esmaily2024new}.
The frequency solver's requirements are identical to a conventional CFD solver, given boundary conditions and a volumetric mesh.
Our results showed that the frequency solver provided accurate outlet flow rate and pressure results within 5\% of a conventional CFD solver while significantly reducing simulation time.
However, the frequency solver's simulation cost scales quadratically with the number of simulated frequency modes due to the nonlinearity of the Navier-Stokes equations.
This cost scaling diminishes the cost-saving benefits of the frequency solver when a large number of frequency modes is needed for a simulation. 

In this study, we will use the harmonic balance method to resolve the cost-scaling issue of the frequency solver.
The harmonic balance method is a frequency domain method used to calculate the response of an unsteady nonlinear system without the need for time integration.~\cite{ROUPHAEL2014179}
It is most widely used for simulating strong nonlinear behaviors such as mechanical vibration and nonlinear electrical circuits~\cite{WANG202327,SUAREZ201665}. 
For fluid dynamics applications, harmonic balance CFD simulations have been utilized in aerospace and turbomachinery studies using the finite volume and the finite difference methods~\cite{mcmullen2001acceleration, jameson2002application, gopinath2007three, sicot2012time,hall2013harmonic,arbenz2017comparison, hupp2016parallel}.
The use of the harmonic balance method for simulating cardiovascular flow has not been thoroughly explored~\cite{koltukluouglu2019harmonic}.
In this paper, we use the harmonic balance form of the Navier-Stokes equations as proposed in a previous paper~\cite{hall2002computation}. 
This form of the equations is a direct transform of the frequency Navier-stokes equations used in our previous study~\cite{esmaily2024new}.
The harmonic balance form has the advantage of using only real-valued variables and allows for approximately linear cost-scaling with the number of frequency modes.
Solving this form of equations results in considerably lower computation costs when a large number of frequency modes is simulated. 

This paper represents the first instance of combining the harmonic balance method with the stabilized finite element method to solve the Navier-Stokes equations. 
Using the proposed harmonic balance finite element solver, we will simulate three different cardiovascular cases with patient-specific anatomy: a Glenn procedure pulmonary flow, a cerebral flow, and a coronary artery flow.
Subsequently, we will compare the results with those obtained from a conventional time-stepping CFD solver. 
These cases hold unique clinical significance and hemodynamic flow patterns and will serve as compelling examples to showcase the effectiveness of the proposed harmonic balance solver. 
We aim to demonstrate that the proposed solver is ready for application in various physiological cases, thereby paving the way for its potential use in future studies.
The paper is organized as follows: we will first introduce the derivation of the harmonics-balance Navier-Stokes equations and present the stabilized finite element method. 
Next, we will provide comprehensive details about the solution procedure to solve the system of equations. 
We will present results obtained from the three cardiovascular cases and use them to demonstrate the harmonic balance solver's cost and accuracy compared to the conventional solver. 
Finally, we will discuss the implications of the results, the existing limitations, and potential future work.

\section{Methods}

This section will outline the numerical methods used to create the harmonic balance CFD solver.  We will present the derivation of the harmonic balance form of the Navier-Stokes equations. 
We will then present the stabilized finite element weak form utilized in this study and briefly discuss other stabilization methods. 
Finally, we will present the solution procedure for solving the resulting system of equations to achieve optimal speed and cost scaling. 
This section will provide sufficient details for readers to adapt an existing conventional finite element solver to solve the harmonic balance Navier-Stokes equations.

\subsection{Harmonic balance Navier-Stokes equations}
In this study, we will use the incompressible flow assumption since we are solving blood flow. 
We also consider blood as Newtonian since the level of shear and size of the major blood vessels is significantly larger than those encountered in capillaries. 
We are also only considering stationary walls. 
Further research is needed to apply this method to solve fluid-structure interactions.

With these assumptions, the three-dimensional incompressible Navier-Stokes equations in the Cartesian coordinate system and the time domain is formulated as  
\begin{equation}
    \begin{alignedat}{3}
    \rho \frac{\partial \hat u_i}{\partial t} + \rho \hat u_j \frac{\partial \hat u_i}{\partial x_j} &= - \frac{\partial \hat p}{\partial x_i} + \mu \frac{\partial^2 \hat u_i}{\partial x_j \partial x_j}, \hspace{0.5in} && \text{in}\;\;\; && \Omega\times(0,N_cT] \\
    \frac{\partial \hat u_i}{\partial x_i}  &= 0,\;\;\; && \text{in} && \Omega\times(0,N_cT]  \\
    \hat u_i &= \hat g_i,  && \text{on} && \Gg\times(0,N_cT]  \\
    -\hat pn_i + \mu \frac{\partial \hat u_i}{\partial x_j}n_j  &= \hat h n_i, && \text{on} && \Gh\times(0,N_cT] 
    \end{alignedat}
 \label{NS-time}
\end{equation}
where subscripts $i$ and $j$ denotes the spatial orientation where $i,j = 1,2,3$. $\rho \in \mathbb R^+$ and $\mu \in \mathbb R^+$ are the fluid density and dynamic viscosity, respectively, and $\hat u_i(\bl x, t) \in \mathbb R$ and $\hat p(\bl x, t) \in \mathbb R$ are the unknown fluid velocity and pressure, respectively, at position $\bl x$ and time $t$. 
$\Omega$ denotes the entire computational domain, and $\Gh$ and $\Gg$ are the portion of the boundary $\Gamma = \partial \Omega = \Gg \bigcup \Gh$ where a Neumann and Dirichlet boundary condition is imposed, respectively.
$N_c$ and $T$ are the number and duration of cardiac cycles. 
$\hat g_i(\bl x, t) \in \mathbb R$ and $\hat h(\bl x, t) \in \mathbb R$ are the imposed Dirichlet and Neumann boundary conditions, respectively, and $n_j$ denotes the outward normal vector to $\Gamma$. 
Note that the Neumann boundary condition is simplified by decoupling the velocity components, which is common for cardiovascular applications~\cite{Vignon2010625, Bazilevs20093534,esmaily2024new}.

To write \eqref{NS-time} in the frequency domain, we take advantage of the underlying frequency of the flow field, $\omega  = 2\pi/T$, calculated from the duration of the cardiac cycle. 
With this knowledge, we can write the solution and boundary conditions in a discretized and truncated frequency form up to the $M$-th frequency mode as
\begin{equation}
\begin{split}
    \hat u_i(\bl x,t) &= \sum_{|m|<M} u^*_{mi}(\bl x) e^{\hat im\omega t}, \\ 
    \hat p(\bl x,t) &= \sum_{|m|<M} p^*_m(\bl x) e^{\hat im\omega t},      \\ 
    \hat g_i(\bl x,t) &= \sum_{|m|<M} g^*_{mi}(\bl x) e^{\hat im\omega t},  \\ 
    \hat h(\bl x,t)& = \sum_{|m|<M} h^*_m(\bl x) e^{\hat im\omega t},    
\end{split}
    \label{upgh_def}
\end{equation} 
where
\begin{align*}
     g^*_{mi}(\bl x) & \coloneqq \frac{1}{T}\int_0^T \hat g_i(\bl x,t) e^{-\hat im\omega t} dt,  \\ 
     h^*_m(\bl x) & \coloneqq \frac{1}{T}\int_0^T \hat h(\bl x, t) e^{-\hat im\omega t} dt,
\end{align*} 
are the Fourier coefficients computed such that the spectral content of the boundary conditions is exact up to the $M$-th mode. 
This way, the error associated with discretizing the boundary conditions in \eqref{upgh_def} scales with the frequency content of the truncated terms and is minimized. The asterisk ($^*$) denotes the variables in frequency and $\hat i = \sqrt{-1}$.

To write the frequency Navier-Stokes equations in a compact matrix form, we organize unknowns and boundary conditions in vectors as 
\begin{equation}
    \bl u^*_i(\bl x) \coloneqq \left[ 
    \begin{matrix} u^*_{(-M+1)i} \\ \vdots\\ u^*_{0i} \\ \vdots \\ u^*_{(M-1)i}\end{matrix} \right], \hspace{0.15in}
    \bl p^*(\bl x) \coloneqq \left[ \begin{matrix} p^*_{-M+1} \\ \vdots \\ p^*_{0} \\ \vdots \\ p^*_{M-1} \end{matrix} \right],  \hspace{0.15in}
    \bl g^*_i(\bl x) \coloneqq \left[ \begin{matrix} g^*_{(-M+1)i}\\ \vdots\\ g^*_{0i}\\ \vdots\\ g^*_{(M-1)i} \end{matrix} \right], \hspace{0.15in}
    \bl h^*(\bl x) \coloneqq \left[ \begin{matrix} h^*_{-M+1}\\ \vdots \\ h^*_{0} \\ \vdots \\ h^*_{M-1} \end{matrix} \right]. 
\label{upghDef}
\end{equation} 

With these definitions, \eqref{NS-time} is expressed in the frequency domain as
\begin{equation}
\begin{alignedat}{3}
    \rho \bl \Omega \bl u^*_i + \rho \bl A_j \frac{\partial \bl u^*_i}{\partial x_j} &= - \frac{\partial \bl p^*}{\partial x_i} + \mu \frac{\partial^2 \bl u^*_i}{\partial x_j \partial x_j},\hspace{0.5in} && \text{in}\;\;\; && \Omega  \\
    \frac{\partial \bl u^*_i}{\partial x_i}  &= 0, && \text{in} && \Omega \\
    \bl u^*_i &= \bl g^*_i,  && \text{on} && \Gg  \\
    -\bl p^*n_i + \mu \frac{\partial \bl u^*_i}{\partial x_j}n_j  &= \bl h^* n_i, && \text{on} && \Gh   
\end{alignedat}
\label{NS-spec}
\end{equation}
where the frequency-source matrix, $\bl \Omega$, is a diagonal matrix with imaginary entries, $\bl \Omega(j,j) \coloneqq \hat i(j-M)\omega$, and $\bl A_i$ is a convolution matrix given as $\bl A_i(j,k) \coloneqq u^*_{(j-k)i}$, $|j-k|<M$. The constraint $|j-k|<M$ results in zero entries on the top right and bottom left of $\bl A_i$, which eliminates the aliasing effect by preventing the generation of frequency that is equal or higher than $M$ through the nonlinear convective acceleration term.

In our previous paper, we used the finite element method with a specialized stabilization scheme to solve equations~\eqref{NS-spec}~\cite{esmaily2024new}.
The computation of the Galerkin convective acceleration term as well as the specialized stabilization terms requires $O(M^2)$ operations, owing to the full convective mode-coupling matrix $\bl A$. 
As a result, the simulation cost scales quadratically with the number of frequency modes. 
This second-order scaling diminished the frequency formulation's cost advantage as the number of frequencies captured in a simulation increased. 

In order to achieve first-order cost scaling, we use a transformation of the frequency Navier-Stokes equations proposed in a previous harmonic balance method paper~\cite{hall2002computation}. We write the vectors containing the variables in frequency as discrete Fourier operations
\begin{equation}
    \bl u^*_i(\bl x) \coloneqq \bl E \bl u_i(\bl x), \hspace{0.15in}
    \bl p^*(\bl x) \coloneqq \bl E \bl p(\bl x),  \hspace{0.15in}
    \bl g^*_i(\bl x) \coloneqq \bl E \bl g_i(\bl x), \hspace{0.15in}
    \bl h^*(\bl x) \coloneqq \bl E \bl h(\bl x),
\label{upghdft}
\end{equation} 
where $\bl E$ is the discrete Fourier transform matrix defined as $\bl E(j,k) \coloneqq (1/N) e^{-2\pi\hat i jk/N}$, where $j, k = 0,\dots,N-1$ and $N = 2M-1$.
The unknowns and boundary conditions are transformed to vectors that contain variables at $N$ equally spaced time points in a cycle, given as
\begin{equation}
    \bl u_i(\bl x) \coloneqq \left[ 
    \begin{matrix} u_{0i} \\ u_{1i} \\ \vdots\\ u_{(N-1)i}\end{matrix} \right], \hspace{0.15in}
    \bl p(\bl x) \coloneqq \left[ \begin{matrix} p_{0} \\ p_{1} \\ \vdots \\ p_{(N-1)} \end{matrix} \right],  \hspace{0.15in}
    \bl g_i(\bl x) \coloneqq \left[ \begin{matrix} g_{0i}\\ g_{1i} \\ \vdots\\ g_{(N-1)i} \end{matrix} \right], \hspace{0.15in}
    \bl h(\bl x) \coloneqq \left[ \begin{matrix} h_{0}\\ h_1 \\ \vdots \\ h_{(N-1)} \end{matrix} \right].
\label{upghharm}
\end{equation} 
Plugging equations~\eqref{upghdft} in~\eqref{NS-spec} and multiplying all equations by the discrete inverse Fourier transform matrix, $\bl E^{-1}$, gives us the final form of the harmonic balance Navier-Stokes equations as

\begin{equation}
\begin{alignedat}{3}
    \rho \bl H \bl u_i + \rho \bl U_j \frac{\partial \bl u_i}{\partial x_j} &= - \frac{\partial \bl p}{\partial x_i} + \mu \frac{\partial^2 \bl u_i}{\partial x_j \partial x_j},\hspace{0.5in} && \text{in}\;\;\; && \Omega  \\
    \frac{\partial \bl u_i}{\partial x_i}  &= 0, && \text{in} && \Omega \\
    \bl u_i &= \bl g_i,  && \text{on} && \Gg  \\
    -\bl p n_i + \mu \frac{\partial \bl u_i}{\partial x_j}n_j  &= \bl h n_i, && \text{on} && \Gh   
\end{alignedat}
\label{NS-harm}
\end{equation}
where $\bl H$ is a zero diagonal time point-coupling matrix where $\bl H = \bl E^{-1}\bl\Omega\bl E$, and $\bl U_i$ is a diagonal convective acceleration matrix where $\bl U_i = \mathrm{diag}\left(\bl u_i \right)$.

\paragraph{Remarks:} 
\begin{enumerate}
    \item Equations \eqref{NS-harm} can be directly obtained from~\eqref{NS-time} if the Eulerian acceleration term is discretized using spectral basis \cite{moin2010fundamentals}. This discretization can be differentiated from the conventional time stepping approach that uses polynomial basis with some form of finite differencing to express the time derivative in a discrete setting. 
    \item The harmonic balance equations in~\eqref{NS-harm} is a direct transformation of the frequency equations in~\eqref{NS-spec}, using a discrete inverse Fourier transform operator. The solutions retain spectral accuracy as a result of the truncated frequency series used to represent the unknowns. 
    \item The discretization of the acceleration term arises from expressing the solution and boundary conditions in a truncated frequency form. This approach assumes the solution consists solely of integer multiples of a fundamental frequency, $\omega$. This assumption is reasonable for cardiovascular flows because the energy in each mode is either directly injected by the boundary condition or produced through interaction between modes. However, in some instances, another distinct frequency, $\omega_g$, may arise from geometry-induced instabilities or turbulence, which is not integer multiples of $\omega$. This typically occurs when the Reynolds number in the fluid domain reaches or exceeds the transient flow limit, typically between $2,000$ and $3,000$. In these cases, the solutions will fail to converge to steady values. Exploration of such unstable flow scenarios are left for the future studies.  
    \item Apart from the time-coupling term, $\rho \bl H \bl u_i$, in equations~\eqref{NS-harm}, all other terms can be calculated at each time point independently of other time points. This allows for the existing code structure used to solve conventional time-domain Navier-Stokes equations to remain largely unchanged. Additionally, since $\bl H$ is formulated with two discrete Fourier transform matrices, a fast Fourier transform routine can be employed, resulting in a cost scaling of $O(N\log{(N)})$. The implementation details will be discussed in the following sections. 
\end{enumerate}

\subsection{Stabilized finite element method}

The discrete finite element form of the Navier-Stokes equations, namely the Galerkin form, needs to be stabilized to mitigate nonphysical solution oscillation encountered in convection-dominant flows~\cite{hughes1987recent,codina2000stabilization, hughes2010stabilized}. In addition, a pressure stabilization term is required to allow the use of equal order shape functions for velocity and pressure~\cite{hughes1979multidimentional, brooks1982streamline}. 
While various stabilization methods have been proposed to counteract these difficulties in the past~\cite{tezduyar2000finite, hsu2010improving,takizawa2018stabilization,wang2023photogrammetry}, we chose Galerkin/least-squares stabilization (GLS) in this study, which successfully counteract those issues while producing a discrete form that remains compact. 
Additionally, various properties of this GLS formulation can be established formally to ensure its convergence rate and numerical stability~\cite{esmaily2024new,esmaily2025augmented}. 
We make the use the same stabilization approach for the conventional time formulation to allow for an apple-to-apple comparison of the two formulations in this study. 

The GLS method for the harmonic balance Navier-Stokes equations is stated as finding $\bl u_i^h$ and $\bl p^h$ such that for any test functions $\bl w_i^h$ and $\bl q^h$, we have
\begin{equation}
\begin{split}
       & \overbrace{\inner{\bl w^h_i,\rho \bl H \bl u^h_i + \rho \bl U_j \frac{\partial \bl u^h_i}{\partial x_j}}_\Omega + \inner{\frac{\partial\bl w^h_i}{\partial x_j}, -\bl p^h \delta_{ij} + \mu \frac{\partial \bl u^h_i}{\partial x_j}}_\Omega  + \inner{\bl q^h, \frac{\partial \bl u^h_i}{\partial x_i}}_\Omega}^\text{Galerkin's terms} \\
       + &\underbrace{\inner{\res_i(\bl w^h_i, \bl q^h),\frac{\bl \tau}{\rho} \res_i(\bl u^h_i, \bl p^h)}_\Ot}_\text{Least-square stabilization terms} = \underbrace{ \inner{\bl w_i^h,  \bl h n_i }_\Gh,}_\text{Neumann BCs. terms}
\end{split}
\label{NS-weak}
\end{equation}
where  
\begin{equation*}
    \res_i(\bl u^h_i, \bl p^h) \coloneqq \rho \bl H \bl u_i^h + \rho \bl U_j \frac{\partial \bl u^h_i}{\partial x_j} + \frac{\partial \bl p^h}{\partial x_i} - \mu \frac{\partial^2 \bl u^h_i}{\partial x_j \partial x_j},
\end{equation*}
is the momentum equation differential operator. 
The inner product notation for two vector functions $\bl f(\bl x)$ and $\bl g(\bl x)$ over $S$ is defined as
\begin{equation}
    \inner{\bl f, \bl g}_S \coloneqq \int_S \bl f^{\intercal} \bl g dS.
    \label{inner_def}
\end{equation} 
$\Ot$ on the stabilization terms indicates elemental integration. 
Note that the terms corresponding to the Dirichlet boundary conditions are not included in \eqref{NS-weak} as $\bl u_i^h=\bl g_i$ and $\bl w_i^h=\bl 0$ on $\Gg$ are directly built into the solution and trial functions spaces, respectively. 

The stabilization parameter $\bl \tau$ in \eqref{NS-weak} plays a vital role in the performance of the GLS method. It is defined such that in the steady limit when $\omega=0$, the equations recover the traditional GLS stabilization parameter as in the conventional finite element formulations~\cite{hughes1986generalized,shakib1991new}. 
Since the convective matrix $\bl U_j$ is diagonal, $\bl \tau$ becomes a diagonal matrix, $\bl \tau = \rm{diag}\left(\tau_0, \tau_1, \dots, \tau_{(\textit{N}-1)}\right)$, which entries can be calculated pointwise for each time point as
\begin{equation}
    \tau_n = \left( u^h_{ni} \cdot \bl {\xi}  u^h_{ni}+ C_{\mathrm{I}}\nu^2\bl {\xi}:\bl {\xi}\right)^{-\frac{1}{2}}, n\in[0, N)
    \label{tau}
\end{equation}
where $\nu$ is the kinematic viscosity, $\boldsymbol{\xi}$ is the covariant tensor obtained from the mapping of the physical-parent elements, and $C_{\mathrm{I}}$ is a shape-function-dependent constant, which is 3 in our study that utilizes tetrahedral elements to mesh all geometries. 

\paragraph{Remarks:} 
\begin{enumerate}
    \item The least-square stabilization terms encapsulate two traditional stabilization terms, the streamline upwind Petrov/Galerkin (SUPG) term $\rho \bl U_j\frac{\partial \bl w^h_i}{\partial x_j} \bl \tau \bl U_k\frac{\partial \bl u^h_i}{\partial x_k}$, and the pressure stabilized Petrov/Galerkin (PSPG) term $\frac{1}{\rho} \frac{\partial \bl q^h}{\partial x_i} \bl \tau \frac{\partial \bl p^h}{\partial x_i}$. However, relying solely on the SUPG/PSPG terms will lead to solution instabilities. Including the full differential operator $\res$, which contains the $\bl H$ terms, is essential for ensuring the stability of this method. This approach contrasts with the stabilization for traditional time-stepping solvers, where the acceleration contribution in $\res$ is often dropped from the formulation. 
    \item Expanding the least-square stabilization terms reveals that this stabilization method resembles the residual-based variational multiscale method (RBVMS)~\cite{bazilevs2007variational,akkerman2008role,bazilevs2009computational}, where the subscale velocity $\bl u_i'$ is included only on the time-coupling term and the velocity gradient in the convection term, i.e.
 $$\inner{\bl w^h_i,\rho \bl H \left(\bl u^h_i + \bl u_i'\right) + \rho \bl U_j \frac{\partial \left(\bl u^h_i + \bl u_i'\right)}{\partial x_j}}_\Omega.$$ 
    Note that the diffusive term naturally drops from $\res$ when linear interpolation functions are used.
    We also tested incorporating the full RBVMS terms, similar to the stabilization method in~\cite{bazilevs2007variational}, by adding $\bl u_i'$ to all velocity terms. However, the solution convergence rate achieved with RBVMS is slightly faster than that with GLS. The overall simulation cost is comparable since RBVMS formulation contains more terms which slightly increased the matrix-vector multiplication cost. 
    \item The GLS method yields accurate results when the element Womersley number $\beta \coloneqq h\sqrt{\frac{(N-1)\omega}{2\nu}}$ is small, which is valid in cardiovascular simulations where the cardiac cycle frequency is small and the mesh resolution is sufficiently large to capture the traveling wave produced by unsteady nature of the flow. In a prior study we have shown that the performance of the GLS method in high-Womersley flow regime can be slightly improved if a two-parameter specialized stabilization method is adopted~\cite{esmaily2025augmented}. 
    \item The conventional design of the stabilization parameter $\bl \tau$ for unsteady problems includes a time-discretization term that in the limit of $\Delta t \to 0$ makes it proportional to $\Delta t$. In conventional methods this term is added to permit solution convergence when a small time step size is used. The downside of adding this term to $\tau$, however, is that it produces a time-inconsistent method~\cite{jia2023time} that fails to converge as $\Delta t \to 0$. The definition in \eqref{tau} eliminates this $\Delta t$-dependent term, hence producing a consistent technique. The reason we are able to remove that term and not encounter convergence issues is that our equivalent $\Delta t$, namely $2\pi/(N\omega)$, is large enough (owing to the use of a few modes for time discretization) that $(1/\Delta t)^2$ can be neglected in comparison to the convective and diffusive terms appearing in \eqref{tau}. Note that the relative importance of those terms can be measured by the element Womersley number $\beta$, which for the representative cases considered in this study remains small $\beta < 0.2$. Our numerical experiments confirm this argument as no difference was observed in the convergence or the final solution when $(1/\Delta t)^2$ was included under root square in \eqref{tau}.
\end{enumerate}

\subsection{Solution Procedure}
This section provides an overview of the solution procedure used to solve equation~\eqref{NS-weak}. We assume all interpolation functions are linear in the following derivations, consistent with what is commonly used in finite element solvers specialized for cardiovascular flow. As a result, the second derivative terms that ought to computed at the Gauss quadrature points are dropped from the formulation. 

It is crucial to note that all the unknowns for all time points are interdependent and must solved simultaneously in our harmonic balance formulation. 
That results in $4N$ unknowns to represent velocity and pressure on each interior node. 

The Newton-Raphson method is employed for iterative root-finding. Additionally, we apply matrix splitting to the tangent matrix, enabling the use of fast Fourier transforms during matrix-vector multiplications. A pseudo time-stepping scheme, developed in our previous study, is utilized here to improve the convergence of the solver~\cite{esmaily2024new}.  

The test functions $\bl w^h_i$ and $\bl q^h$ in \eqref{NS-weak} are discretized in space using
\begin{align*}
    \bl w^h_i(\bl x) & = \sum_A \bl c^{\rm w}_{Ai} N_A(\bl x), \\ 
    \bl q^h(\bl x) & = \sum_A \bl c^{\rm q}_{A} N_A(\bl x), 
\end{align*}
where $N_A(\bl x)$ is the interpolation function associated with node $A$. 
Since $\bl w^h_i$ and $\bl q^h$ are arbitrary functions, \eqref{NS-weak} must hold for any $\bl c^{\rm w}_{Ai}$ and $\bl c^{\rm q}_A$. 
That permits us to obtain a system of equations from \eqref{NS-weak}, which are
\begin{equation}
\begin{split}
      \bl r^{\rm m}_{Ai} = & \int_\Omega \left(\rho N_A \bl H \bl u^h_i + \rho N_A \bl U_j \frac{\partial \bl u^h_i}{\partial x_j} - \frac{\partial N_A}{\partial x_i}\bl p^h + \mu \frac{\partial N_A}{\partial x_j} \frac{\partial \bl u^h_i}{\partial x_j} \right)d\Omega  - \int_\Gh N_A\bl h n_i d\Gamma  \\
       + & \int_\Ot\left(- \bl H N_A + \bl U_j \frac{\partial N_A}{\partial x_j}\right)  \bl \tau \res_i(\bl u^h_i, \bl p^h) d\Omega = \bl 0, \hspace{1.3cm} A\in\eta-\eta_g,\; i=\text{1, 2, and 3,}  \\
       \bl r^{\rm c}_A = & \int_\Omega N_A \frac{\partial \bl u^h_i}{\partial x_i}d\Omega + \int_\Ot \frac{1}{\rho}\frac{\partial N_A}{\partial x_i}\bl \tau \res_i(\bl u^h_i, \bl p^h) d\Omega = \bl 0, \hspace{1cm} A\in\eta. 
\end{split}
       \label{Rmc}
\end{equation}
In \eqref{Rmc}, $\eta$ and $\eta_g$ denote the set of nodes in the entire domain $\Omega$ and those located on the Dirichlet boundaries $\Gg$, respectively. It is important to note that the negative sign in front of the $\bl H N_A$ in the least-square term arises from taking the transpose on the left differential operator. Since $\bl H$ is a skew-symmetric matrix and $\bl U_j$ is a diagonal matrix, we can write the transpose as $\bl H^{\intercal} = - \bl H$ and $\bl U_j^{\intercal} = \bl U_j$.

The velocity and pressure are discretized using the same interpolation functions as those of the test functions, namely
\begin{align*}
    \bl u_i^h(\bl x) & = \sum_A \bl u^{\rm d}_{Ai} N_A(\bl x), \\ 
    \bl p^h(\bl x) & = \sum_A \bl p^{\rm d}_{A} N_A(\bl x), 
\end{align*}
where $\bl u^{\rm d}_{Ai}$ and $\bl p^{\rm d}_A$ contain velocity and pressure, respectively, at all time points at node $A$.
In our implementation, we build the Dirichlet boundary condition into the unknown vector so that $\bl u^{\rm d}_{Ai} = \bl g_i(\bl x_A)$, where $\bl x_A$ is the position of node $A$.  

At the discrete level, our goal is to find $\bl u^{\rm d}_{Ai}$ and $\bl p^{\rm d}_A$ such that all equations in \eqref{Rmc} are satisfied. 
Given that these equations are nonlinear, we accomplish this through an iterative process using the Newton-Raphson iterations, given as solving
\begin{equation}
    \bl y^{(n+1)} = \bl y^{(n)} - \left(\bl L^{(n)}\right)^{-1}\bl r^{(n)}
   \label{NR}
\end{equation}
at each Newton-Raphson iteration $n$ to update the solution from the last iteration $\bl y^{(n)}$ and compute it at the next iteration $\bl y^{(n+1)}$. $\bl r^{(n)}$ is the residual vector calculated based on the unknowns at the previous iteration $n$ using \eqref{Rmc} and $\bl L^{(n)}$ is the tangent matrix calculated from $\bl y^{(n)}$. The unknowns and residuals are written in vectors as
\begin{equation*}
    \bl y \coloneqq \left[\begin{matrix} \bl u^{\rm d}_1 \\ \bl u^{\rm d}_2 \\ \bl u^{\rm d}_3 \\ \bl p^{\rm d} \end{matrix}\right],  \hspace{0.2in}
    \bl r \coloneqq \left[\begin{matrix} \bl r^{\rm m}_1 \\ \bl r^{\rm m}_2 \\ \bl r^{\rm m}_3 \\ \bl r^{\rm c} \end{matrix}\right],  
\end{equation*}

Before presenting the form of the tangent matrix, it is clear that all terms not containing the matrix $\bl H$ in \eqref{Rmc} can be calculated pointwise. 
This means that the residual at each time point can be calculated using only the unknowns associated with that specific time point, provided that the terms involving the time-coupling acceleration matrix $\bl H$ can be calculated separately. With this understanding, we perform matrix splitting on the tangent matrix $\bl L^{(n)}$ to separate time-coupling terms from the non-coupled terms as
\begin{equation}
    \bl L = \frac{\partial \bl r}{\partial \bl y} = \bl C + \bl P,
    \label{splitK}
\end{equation}
where the superscript $(n)$ is dropped to simplify the notation. The from of the tangent matrix $\bl P$ is essentially identical to that of a conventional steady-state finite element Navier-Stokes solver and is
\begin{equation}
\bl P = \left[\begin{matrix}
       \bl K & \bl 0 & \bl 0 & \bl G_1 \\     
       \bl 0 & \bl K & \bl 0 & \bl G_2 \\     
       \bl 0 & \bl 0 & \bl K & \bl G_3 \\     
       \bl D_1 & \bl D_2 & \bl D_3 & \bl L \\     
       \end{matrix} \right],
       \label{ptangent}
\end{equation}
where
\begin{equation}
\begin{split}
       \bl K_{AB} & = \int_\Omega\left(\rho N_A \bl U_k \frac{\partial N_B}{\partial x_k} + \mu \frac{\partial N_A}{\partial x_k} \frac{\partial N_B}{\partial x_k} \bl I \right)d\Omega + \int_\Ot \bl U_k \frac{\partial N_A}{\partial x_k}\rho \bl \tau \bl U_k  \frac{\partial N_B}{\partial x_k} d\Omega,  \\ 
       (\bl G_{AB})_i & = -\int_\Omega \frac{\partial N_A}{\partial x_i} N_B\bl I d\Omega + \int_\Ot \bl U_k \frac{\partial N_A}{\partial x_k} \bl \tau \frac{\partial N_B}{\partial x_i} d\Omega,  \\ 
       (\bl D_{AB})_j & = \int_\Omega N_A \frac{\partial N_B}{\partial x_j} \bl I d\Omega + \int_\Ot \frac{\partial N_A}{\partial x_j}\bl \tau \bl U_k \frac{\partial N_B}{\partial x_k} d\Omega,  \\
       \bl L_{AB} & = \int_\Ot  \frac{1}{\rho} \frac{\partial N_A}{\partial x_k}\bl \tau\frac{\partial N_B}{\partial x_k} d\Omega.  
\end{split}
\label{K_mat}
\end{equation}
It is important to note that all blocks in matrix $\bl P$ are diagonal matrices. Therefore, the computational cost associated with matrix multiplications of $\bl P$ is $O(N)$ using a sparse matrix implementation.

The time point-coupling tangent matrix $\bl C$ handles all terms involving the $\bl H$ matrix. It includes one contribution from the Galerkin's terms and four contributions from the least-squares stabilization terms in \eqref{NS-harm}. Our numerical experiments showed that dropping the stabilization terms' contribution to $\bl C$ has a negligible effect on solution convergence. Therefore, we neglect the tangent contributions from these terms, which simplifies the number of operations needed to compute $\bl C$. It is important to note that this simplification applies only to the calculation of the tangent matrix and not to the residual, ensuring no impact on the final converged solution. The simplified time-coupling tangent matrix takes the following form
\begin{equation}
\bl C = \left[\begin{matrix}
       \bl F & \bl 0 & \bl 0 & \bl 0 \\     
       \bl 0 & \bl F & \bl 0 & \bl 0 \\     
       \bl 0 & \bl 0 & \bl F & \bl 0 \\     
       \bl 0 & \bl 0 & \bl 0 & \bl 0 \\     
       \end{matrix} \right],
    \label{ctangent}
\end{equation}
where
\begin{equation}
       \bl F_{AB} = \bl H \int_\Omega \rho N_A \bl I N_B d\Omega ,  \\ 
\label{F_mat}
\end{equation}
The matrix $\bl H$ is formulated as $\bl H = \bl E^{-1} \bl \Omega \bl E$. The matrix multiplications involving $\bl E$ and $\bl E^{-1}$ can be efficiently handled using a fast Fourier transform (FFT) and an inverse fast Fourier transform (IFFT), respectively. 
This implementation reduces the computational cost of performing matrix-vector multiplications involving $\bl H$ from $O(N^2)$ to $O(N\log(N))$.
When combined with the cost of computing other terms, the overall cost-scaling for the solver becomes approximately $O(N)$. 
This study used the open-source library FFTW~\cite{frigo2005design} to execute the fast Fourier transforms. 

Before solving the final linear system, a pseudo time-stepping scheme is added to improve linear solver convergence and reduce the overall cost.
Pseudo time-stepping schemes have been used previously in solving harmonic balance equations~\cite{hall2002computation} and the specific method used here was proposed in our previous study~\cite{esmaily2024new}.
The implementation of the pseudo time-stepping involves adding 
\begin{equation}
    \inner{\bl w^h_i,\rho\frac{\partial \bl u^h_i}{\partial \tilde t} }_\Omega, 
    \label{sudot}
\end{equation}
to \eqref{NS-weak}, and integrating $\bl u_i^h$ in pseudo-time $\tilde t$ similar to the conventional time formulation of the Navier-Stokes equations. 
Using the generalized-$\alpha$ method~\cite{jansen2000generalized}, the tangent matrix is supplemented by adding 
\begin{equation}
    \left(N_A,\frac{c_1\rho}{\Delta \tilde t}N_B \bl I \right)_\Omega, 
\end{equation}
to $\bl K_{AB}$ in \eqref{K_mat}, where we used $\rho_{\infty} = 0$ resulting in $c_1 = 1.5$.
Since the pseudo time step term goes to zero as the solution converges to a ``steady'' result, it improves the convergence behavior without affecting the final solution. 

The choice of the pseudo time step size $\Delta \tilde t$ is essential for achieving optimal convergence.
Our numerical experiments show that the system of equations will generally converge when the convective Courant–Friedrichs–Lewy (CFL) number calculated from the pseudo time step size is $O(1)$ or less, defined as
\begin{equation}
    C_{\mathrm{CFL}} = \frac{u_c \Delta \Tilde{t}}{h_c},
\end{equation}
where the characteristic length, $h_c$, is taken as the edge length for a regular tetrahedron whose volume equals the mean elemental volume of the entire domain, and the characteristic convective velocity, $u_c$, is taken as the mean velocity at the inlet boundary at the maximum flow rate for the inlet flow profile.
Therefore, the optimal choice of pseudo time step size can be approximated as
\begin{equation*}
    \Delta \Tilde{t} \sim \frac{h_c}{u_c}.
\end{equation*}

\begin{figure} [ht]
    \centering
    \includegraphics[width=0.6\linewidth]{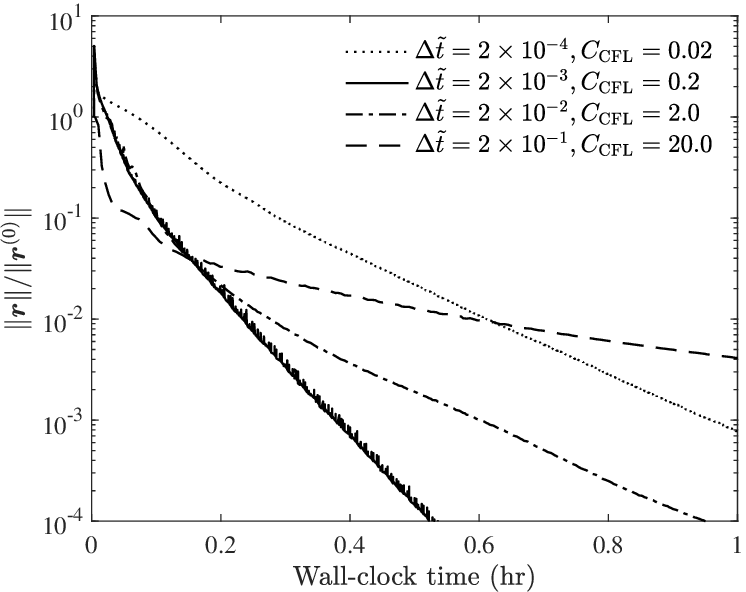}
    \caption{Simulation convergence speed for four pseudo time step sizes, $\Delta \Tilde{t}$. The y-axis is the relative residual compared to the residual at the initialization step, $\|\mbox{\boldmath{$r$}}\|/\|\mbox{\boldmath{$r$}}^{(0)}\|$.}
    \label{fig:sudo}
\end{figure}

To illustrate the impact of the CFL number on convergence, we conducted four identical simulations, varying only the pseudo time step size. 
These simulations are based on the $N=13$ case in Section 3.2.
As shown in Figure~\ref{fig:sudo}, the simulation with a large CFL number converges exceptionally slowly. 
Conversely, an excessively small CFL number also unnecessarily increases simulation costs due to the cost of constructing the linear system at each pseudo step. 
It is important to note that optimal CFL number and pseudo time step size for achieving the fastest convergence may vary depending on other factors, such as the geometry of the case, mesh uniformity, and flow dynamics.
Generally, optimal convergence is attained when the CFL number is around 1. 

A backflow stabilization scheme is used to prevent numerical instabilities caused by reversal flow through $\Gh$, where a Neumann boundary condition is imposed. 
This method is identical to the backflow stabilization scheme used in conventional time-stepping Navier-Stokes finite element solver~\cite{Esmaily2011backflow}. 
The stabilization is achieved by adding
\begin{equation}
    \inner{\bl w^h_i, \frac{\rho}{2}\beta  |\bl U_i n_i|_- \bl u_i^h }_\Gh, 
\end{equation}
to the right-hand side of \eqref{NS-weak}, where 
\begin{equation*}
     |\bl U_i n_i|_- = \frac{\bl U_i n_i - |\bl U_i n_i|}{2},
\end{equation*}
and $\beta \in [0, 1]$ is a user-defined coefficient, which we used $\beta = 0.2$ in all cases.
 
We applied a standard Jacobi preconditioner on the tangent matrix $\bl L$ to improve the conditioning of the linear system~\cite{esmaily2013new}. Since $\bl C$ is a zero diagonal matrix, the preconditioner only consists of the diagonal entries of matrix $\bl P$. To solve the linear system, we utilized the generalized minimal residual method (GMRES)~\cite{saad1986gmres}. Simulations are considered converged once the two-norm of the residual decreases by three orders of magnitude, as established in our precious study~\cite{esmaily2024new}.

\section{Results}
The proposed harmonic balance solver is implemented in our in-house finite-element solver named Multi-physics finite-element solver (MUPFES)~\cite{moghadam2013modular,esmaily2015bi}.
The solver has been verified~\cite{steinman2013variability} and employed for cardiovascular modeling in the past~\cite{esmaily2012optimization,esmaily2015assisted,jia2021efficient}.
This solver is parallelized using a message passing interface (MPI) and a specialized parallelization strategy~\cite{esmaily2013low,esmaily2015impact}.
The workload is parallelized using spatial partitioning by employing ParMETIS library~\cite{METIS}. 
All computations are performed on a cluster of AMD Opteron$^{\rm TM}$ 6378 processors that are interconnected via a QDR Infiniband.
Unless otherwise stated, simulations were run using 288 cores (4.5 nodes) operating at 2.4 GHz. 

In this section, we primarily focus on comparing the accuracy and cost of the proposed harmonic balance CFD solver with the conventional time-stepping CFD solver. 
To ensure a fair comparison, the time-stepping CFD solver also utilizes the GLS stabilization method. However, unlike the proposed solver, the stabilization parameter for the time solver, $\hat \tau$, includes a contribution from the acceleration term to improve numerical stability. 
The exact form was introduced in a prior paper as~\cite{jia2023time}
\begin{equation}
    \hat \tau = \left(\hat \omega^2 + \hat u^h_i \cdot \bl {\xi}  \hat u^h_i + C_{\mathrm{I}}\nu^2\bl {\xi}:\bl {\xi}\right)^{-\frac{1}{2}},
    \label{tauhat}
\end{equation}
where $\hat \omega = \|\frac{\partial{\hat u_i^h}}{\partial t}\|_\Omega/\|\hat u_i^h\|_\Omega$.
The code structure of the time-stepping and harmonic balance solvers are kept as identical as possible for comparison purposes.
The linear solver (GMRES) tolerance is set at $0.03$ for both solvers. 
The time solver Newton-Ralphson iteration tolerance is set at three order of magnitudes, consistent with the convergence criteria of the harmonic balance solver. 

We conducted simulations on three cardiovascular cases using the same mesh for both solvers: a Glenn procedure case involving pulmonary arteries, a cerebral artery case, and a left main coronary artery case. 
These cases were chosen to represent a wide range of scenarios, distinguished both numerically and hemodynamically. 
The geometries and inlet flow rate for these cases were obtained from an online repository~\cite{wilson2013vascular}. 
The Reynolds number for the Glenn case is higher than that of the other two, featuring a smooth inflow profile that can be accurately captured with just a few time points (frequencies).
In contrast, the cerebral arteries case presents a more complex smooth inflow profile with higher geometric complexity. 
Meanwhile, the left main coronary artery case has sharper temporal changes in the inflow profile, which requires a significant number of time points to represent accurately.

The outlets in this study are modeled as Neumann boundaries with fixed physiological pressure values. 
Previous cardiovascular simulations have treated the outlets as resistance boundaries to achieve more realistic outlet responses~\cite{esmaily2012optimization,troianowski2011three}. 
However, implementing the resistance boundary for the harmonic balance solver is beyond the scope of this study.
All walls are assumed to be non-moving and non-slip. 
The whole blood density ($\rho = 1.06$ $\mathrm{g/cm^3}$) and dynamic viscosity ($\mu = 4$ $\mathrm{mPa\cdot s}$) are constant across all cases. 
These results aim to demonstrate the solution validity and cost efficiency of the harmonic balance solver as an alternative to the conventional time-stepping solver. 
Since the geometries and boundary conditions in all cases are physiological, we will also briefly discuss the clinical relevance of each case.

\subsection{Glenn pulmonary flow}

\begin{figure} [h]
    \centering
    \includegraphics[width=0.6\linewidth]{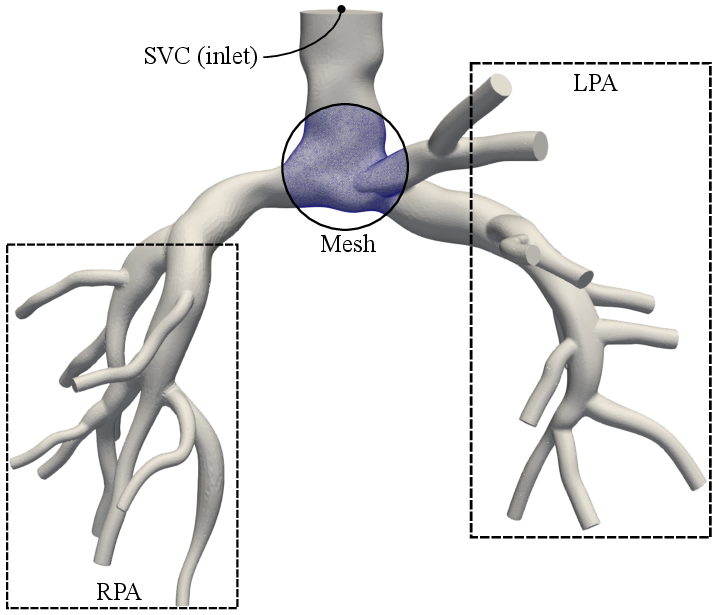}
    \caption{Clinically obtained Glenn pulmonary geometry with a sectional view of the mesh. SVC, superior vena cava; LPA, left pulmonary artery; RPA, right pulmonary artery.}
    \label{fig:pulm_geo}
\end{figure}

The Glenn procedure is an early-stage surgical intervention designed for patients with single ventricle birth defects~\cite{mainwaring1999effect,pridjian1993usefulness}. 
The operation establishes the pulmonary circulation by directly connecting the superior vena cava (SVC) to the pulmonary arteries (PA), as shown in the geometry in Figure~\ref{fig:pulm_geo}.
Prior studies have utilized this geometry to conduct conventional CFD simulations~\cite{troianowski2011three}. 
The inflow boundary condition is taken from the prior study, providing the SVC flow rate as a function of time during a cardiac cycle.
The outlet boundaries are set at a steady pressure of 7.5 mmHg. 
The cardiac cycle duration is 0.6979 seconds, and the peak Reynolds number is around 1,000. 

\begin{figure}[H]
  {\captionsetup{position=bottom,justification=centering}\begin{subfigure}[b]{0.32\textwidth}
    \centering
    \includegraphics[width=\textwidth]{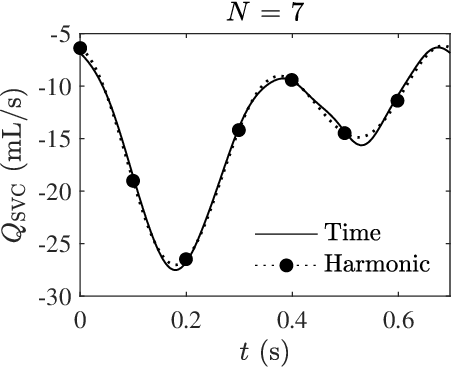}
    \caption[]{{\small}}
    \label{fig:pulm_bc_n7}
  \end{subfigure}
  \hfill
  \begin{subfigure}[b]{0.32\textwidth}
    \centering
    \includegraphics[width=\textwidth]{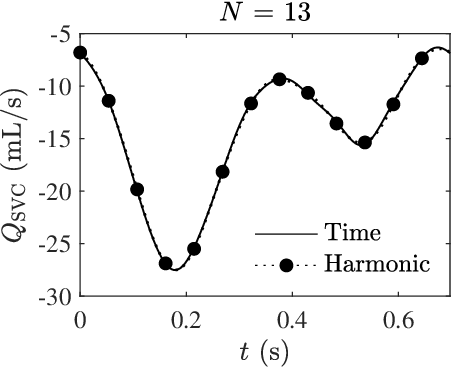}
    \caption[]{{\small}}
    \label{fig:pulm_bc_n13}
  \end{subfigure}
\hfill
    \begin{subfigure}[b]{0.32\textwidth}
    \centering
    \includegraphics[width=\textwidth]{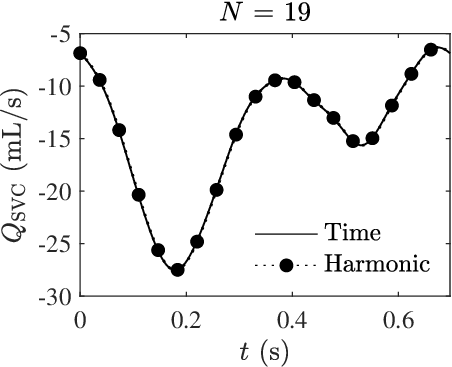}
    \caption[]{{\small}}
  \label{fig:pulm_bc_n19}
  \end{subfigure}}
  \caption{Discretized harmonic balance boundary conditions at each time point (black circle) on the SVC inlet. The dotted line represents the reconstructed profile from the time points. For reference, the solid line is the original flow profile in time. (a) $N = 7$, (b) $N = 13$, (c) $N = 19$.}
  \label{fig:pulm_bc}
\end{figure}

In total, four simulations are performed: one using the conventional time-stepping CFD solver and ran for four cardiac cycles ($4,000$ time steps), and three using the harmonic balance solver with varying numbers of time points ($N = 7$, 13, and 19).
The geometry is discretized with 1,946,676 tetrahedral elements.
For the harmonic balance simulations, the pseudo time step size is $\Delta \Tilde{t} = 5 \times 10^{-3}$ resulting in $C_{\mathrm{CFL}} = 2.5$.
The inlet boundary conditions for the harmonic balance simulations are generated by performing a Fourier transform on the inflow profile, followed by a discrete inverse Fourier transform on the truncated Fourier series, as defined in Equations~\eqref{upgh_def} and~\eqref{upghdft}. 
The resulting inlet boundary conditions are depicted in Figure~\ref{fig:pulm_bc}.
As shown in the boundary condition plots, this inlet profile for this case is smooth and can be accurately captured with a few time points. 
The truncation error for the inflow profile is already less than 1\% for $N = 13$.

\begin{figure} [h]
    \centering
    \includegraphics[width=0.5\linewidth]{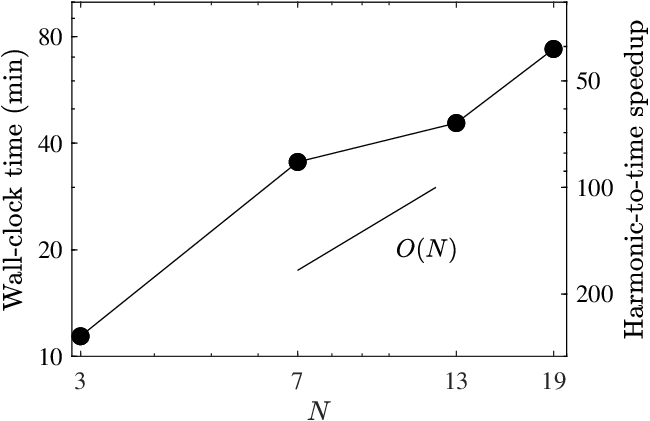}
    \caption{The cost and speed-up of the harmonic balance solver for the pulmonary flow case. The left axis shows the wall-clock time of the harmonic balance simulations. The right axis shows the speedup of the harmonic balance simulations compared to the conventional time-stepping simulation. A linear scaling line is included for reference.}
    \label{fig:pulm_cost}
\end{figure}

Figure~\ref{fig:pulm_cost} illustrates the cost of performing the harmonic balance simulations. The left axis shows the wall-clock time of the simulations, and the right axis indicates the speedup compared to the conventional time-stepping simulation. The conventional simulation took 50 wall-clock hours to complete. 
The figure is presented on a log-log scale to highlight the linear cost scaling of our solver as the number of time points increases. 
The harmonic balance simulation using $N=3$ was included solely to demonstrate cost scaling. 
For the case with the largest number of time points, $N=19$, the harmonic balance simulation is still 40 times faster than the conventional simulation.

\begin{figure}[H]
  {\captionsetup{position=bottom,justification=centering}\begin{subfigure}[b]{0.48\textwidth}
    \centering
    \includegraphics[width=\textwidth]{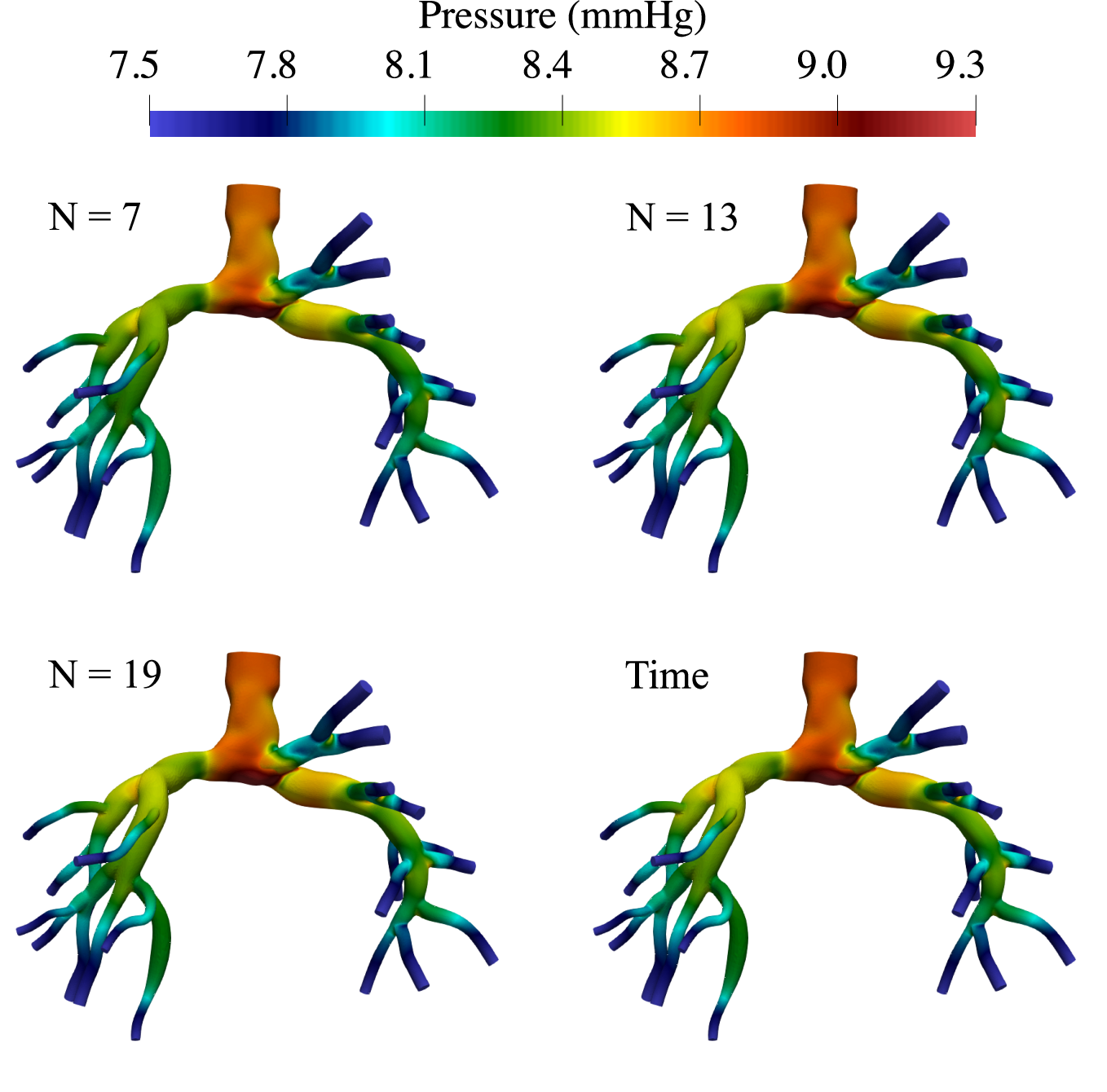}
    \caption[]{{\small}}
    \label{fig:pulm_pres_cont}
  \end{subfigure}
  \hfill
  \begin{subfigure}[b]{0.48\textwidth}
    \centering
    \includegraphics[width=\textwidth]{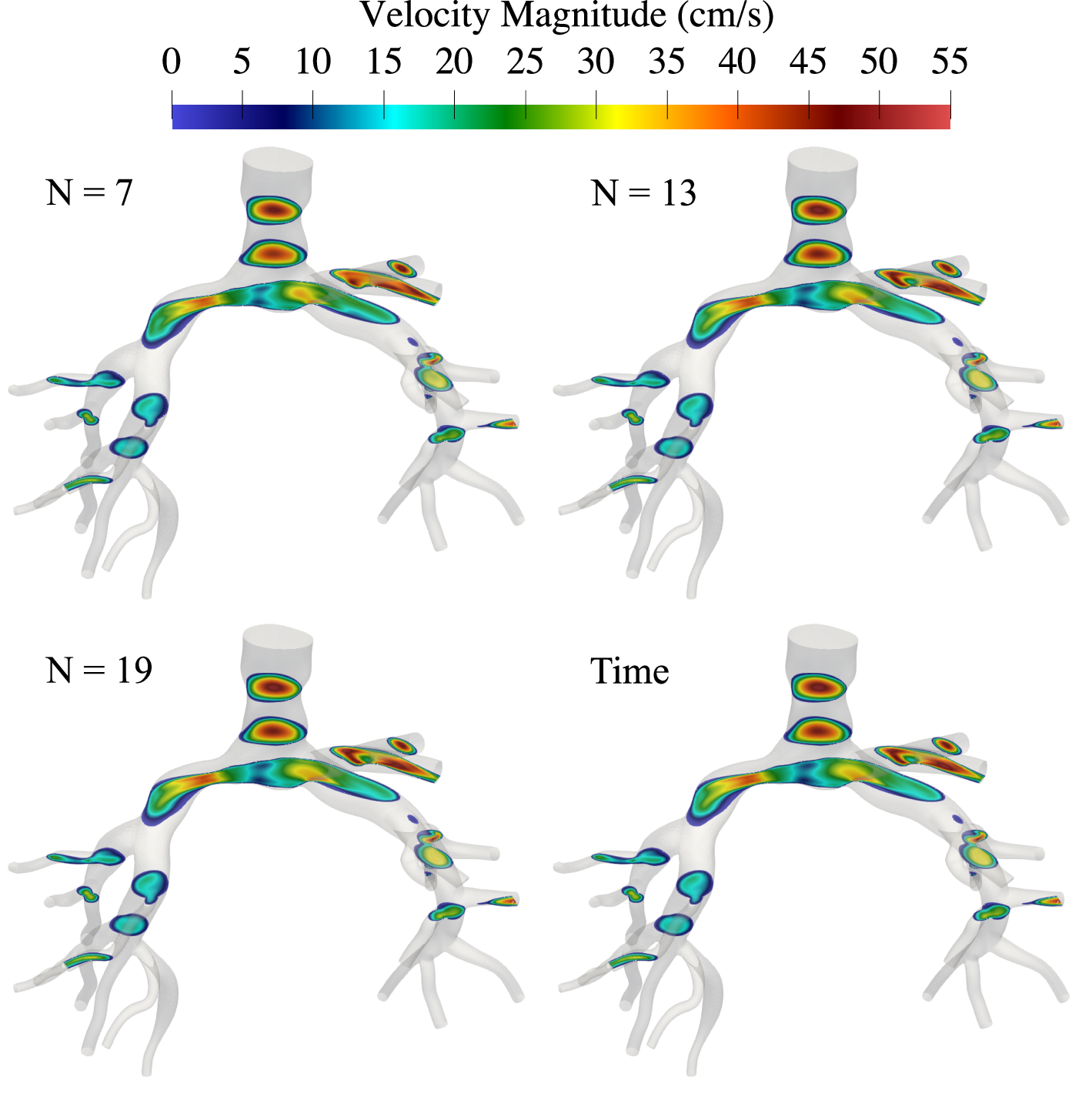}
    \caption[]{{\small}}
    \label{fig:pulm_vel_cont}
  \end{subfigure}}
  \caption{Harmonic balance results using $N = 7$, 13, and 19 compared to the conventional time results for the pulmonary flow case. The (a) pressure and (b) velocity magnitude contours are taken at $t = 0.2$ seconds of the cardiac cycle. }
  \label{fig:pulm_cont}
\end{figure}

To obtain the continuous temporal solution of the flow field, we can perform a Fourier transform on the time-point solutions. We then reconstruct the temporal solution from the obtained Fourier modes, as defined in Equation~\ref{upgh_def}.
Figure~\ref{fig:pulm_cont} illustrates the velocity and pressure fields at $t = 0.2 s$ of the cardiac cycle, when the inlet flow rate reaches its peak.
Qualitative differences are only observed in the $N=7$ case at certain locations in the domain, such as the SVC-PA junction.

In addition to qualitatively examining the velocity and pressure fields, we aim to establish a quantitative measurement for the relative accuracy of the harmonic balance solver compared to the conventional time solver.
To achieve this, we introduce a notation for the integrated relative root mean squared error, $E_X^{f(x)}$, as a measure of result error. 
The error is defined as
\begin{equation}
    E_X^{f(x)} =  \frac{\int_{X}\|f_{\mathrm{harmonic}}(x) - f_{\mathrm{time}}(x)\|dx}{\int_{X}\|f_{\mathrm{time}}(x)\|dx} ,
    \label{eq:re}
\end{equation}
where the dependent variable $f(x)$ can represent either the velocity vector, $\bl u$, or the pressure, $p$. The subscripts indicate whether the results are taken from the harmonic balance solver or the conventional time solver. 
The integral can be performed over the geometry, $\Omega$, the cardiac cycle, $T$, or both. 

\begin{figure}[H]
{\captionsetup{position=bottom,justification=centering}\begin{subfigure}[b]{0.48\textwidth}
    \centering
    \includegraphics[width=\textwidth]{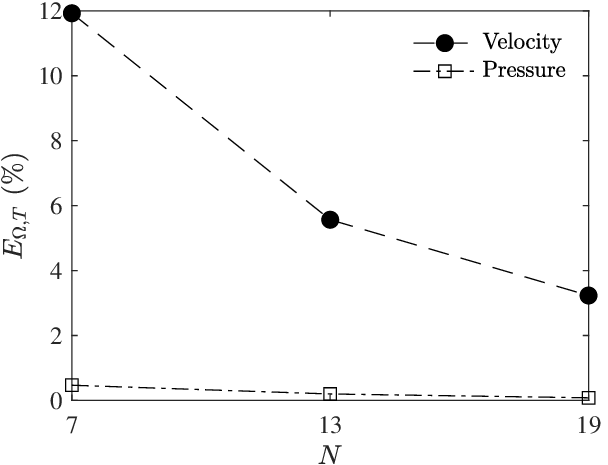}
    \caption[]{{\small}}
    \label{fig:pulm_e_int}
  \end{subfigure}
  \hfill
  \begin{subfigure}[b]{0.48\textwidth}
    \centering
    \includegraphics[width=\textwidth]{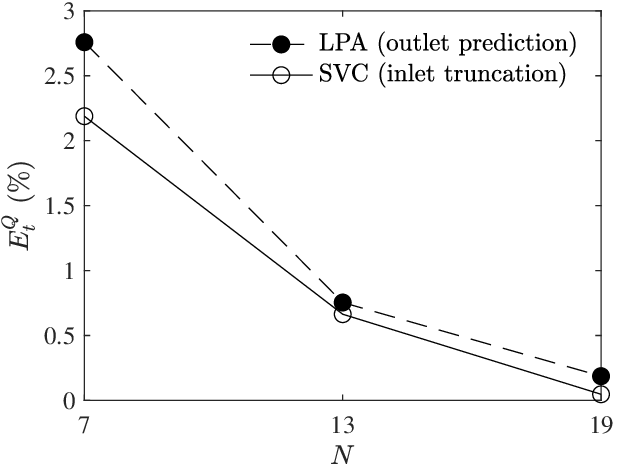}
    \caption[]{{\small}}
    \label{fig:pulm_eQ}
  \end{subfigure}}
  \caption{Pulmonary flow relative root mean square error for (a) the velocity and pressure integrated over the fluid domain and one cardiac cycle and (b) the inlet SVC flow rate (truncation) and the LPA outlet flow rate (prediction). }
  \label{fig:pulm_eint}
\end{figure}

Using this definition, we can calculate the solution error by integrating it over the fluid domain and over one cardiac cycle, and observe the trend as the harmonic balance time point resolution ($N$) increases.
As depicted in Figure~\ref{fig:pulm_e_int}, the error of the variables is integrated both across the fluid domain, $\Omega$, and over one cardiac cycle, $T$. 
The pressure relative error is considerably smaller than the velocity error because the Neumann outlets are set at a fixed value that is much larger than the range of pressure in the solution. 
The velocity error decreases substantially as the number of time points increases, falling to less than $5\%$ when $N = 19$. 
This result also demonstrates that while the inlet flow profile can be well represented with $N = 7$, due to the intermediate Reynolds number of this case ($\mathrm{Re} = 1000$), the flow energy is transferred to the higher frequency modes due to the nonlinear convective effect. 
Therefore, more than seven time points are required to accurately capture the dynamics of the fluid field. 

\begin{figure}[H]
  {\captionsetup{position=bottom,justification=centering}\begin{subfigure}[b]{0.32\textwidth}
    \centering
    \includegraphics[width=\textwidth]{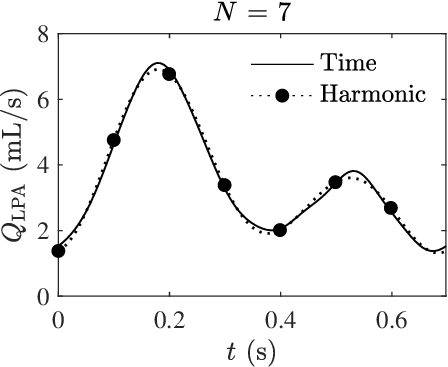}
    \caption[]{{\small}}
    \label{fig:pulm_lpa_n7}
  \end{subfigure}
  \hfill
  \begin{subfigure}[b]{0.32\textwidth}
    \centering
    \includegraphics[width=\textwidth]{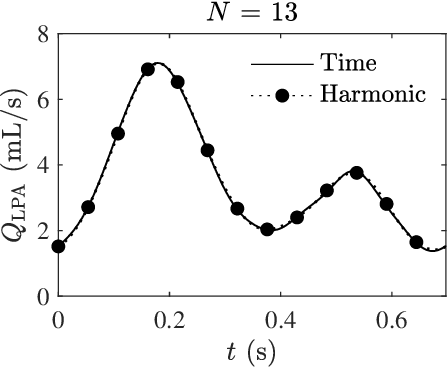}
    \caption[]{{\small}}
    \label{fig:pulm_lpa_n13}
  \end{subfigure}
\hfill
    \begin{subfigure}[b]{0.32\textwidth}
    \centering
    \includegraphics[width=\textwidth]{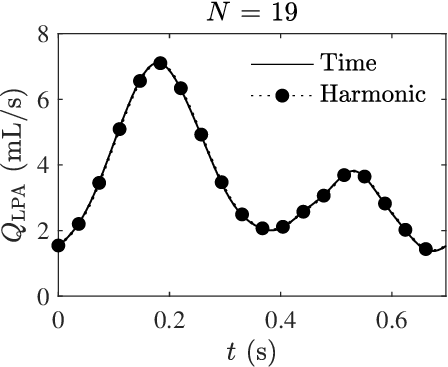}
    \caption[]{{\small}}
  \label{fig:pulm_lpa_n19}
  \end{subfigure}}
  \caption{Left pulmonary artery (LPA) flow rate results obtained with the harmonic balance solver (solid circles) and the with the conventional time solver (solid line). The dotted line represents the reconstructed profile from the harmonic balance time points. (a) $N = 7$, (b) $N = 13$, (c) $N = 19$.}
  \label{fig:pulm_lpa}
\end{figure}

It is important to note that the relative error, $E_{\Omega,T}$, is integrated over both space and time. 
This aggregate nature of the error calculation causes it to serve as an estimate for the entire solution field. 
However, in practice, the parameters of physiological significance are often summarized values at specific locations, such as flow rates at the outlets.
For instance, for patients who underwent the Glenn operation, the flow rate to each side of the pulmonary artery is a critical indicator for potential asymmetric lung development~\cite{cloutier1985abnormal,mendelsohn1994central}. 
Figure~\ref{fig:pulm_lpa} shows the flow rate over one cardiac cycle for the left pulmonary artery (LPA), $Q_{\mathrm{LPA}}$. This value is obtained by integrating the velocity normal to the LPA outlet surface. 
Notably, the flow profiles derived from the harmonic balance results closely resemble those from the conventional time results, even in the simulation using only seven time points ($N = 7$). 
We can compute the relative error of the LPA flow rate and compare it with the inlet SVC flow rate truncation.
As demonstrated in Figure~\ref{fig:pulm_eQ}, the inlet truncation error provides an excellent estimate of the solution error for the LPA flow rate. Specifically, the error is less than $3\%$ when $N = 7$ and drops well below $1\%$ when the number of time points exceeds $N = 13$.

\subsection{Cerebral flow}

\begin{figure} [H]
    \centering
    \includegraphics[width=0.3\linewidth]{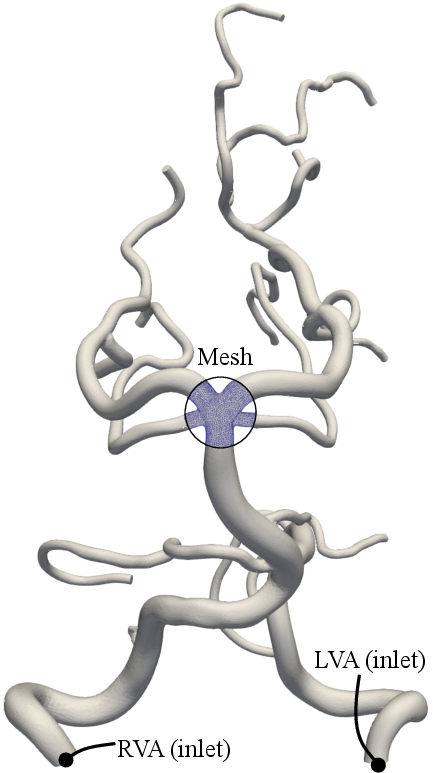}
    \caption{Clinically obtained cerebral arteries geometry used in this study with a sectional view of the mesh resolution. LVA, left vertebral artery; RVA, left vertebral artery.}
    \label{fig:cere_geo}
\end{figure}

The second case presented in the study examines blood flow in the cerebral arteries, which supply approximately 20\% of the brain's blood. 
The geometry includes two vertebral arteries (VA) inlets that converge into the basilar artery, as shown in Figure~\ref{fig:cere_geo}.
This geometry was previously used in conventional computational fluid dynamics (CFD) simulations to analyze the mixing of flows from each side of the vertebral arteries~\cite{bockman2012fluid}. 
The inlet flow rate profile mirrors that of the previous study, with a period of 1 second. 
The outlet pressure is set to the average value reported in the prior study at 115 mmHg~\cite{bockman2012fluid}. 
The peak Reynolds number for this flow is approximately 120.
Unlike the simpler pulmonary flow, the inflow profile for cerebral circulation is more complex and requires additional time points for accurate representation. 
In fact, the Fourier frequency content of the inflow includes exactly 10 unique non-zero values. 
Therefore, the inflow profile can be exactly represented with $N = 19$ numbers of time points. 

\begin{figure}[h]
    \centering
  {\captionsetup{position=bottom,justification=centering}\begin{subfigure}[b]{0.32\textwidth}
    \centering
    \includegraphics[width=\textwidth]{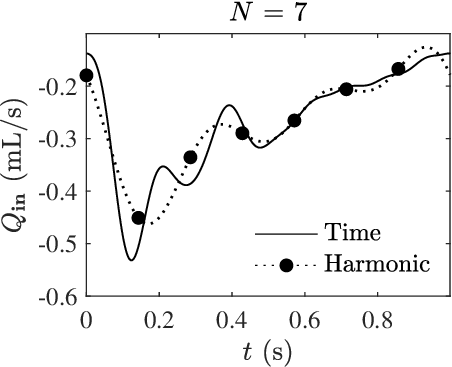}
    \caption[]{{\small}}
    \label{fig:cere_bc_n7}
  \end{subfigure}
\hfill  \begin{subfigure}[b]{0.32\textwidth}
    \centering
    \includegraphics[width=\textwidth]{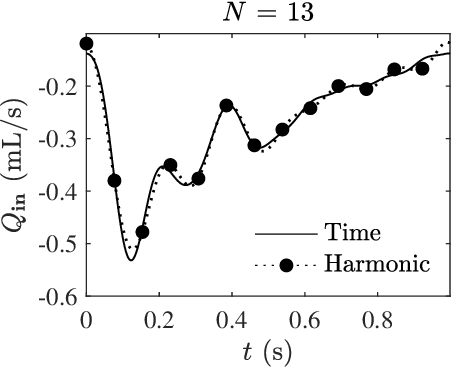}
    \caption[]{{\small}}
    \label{fig:cere_bc_n13}
  \end{subfigure}
\hfill
    \begin{subfigure}[b]{0.32\textwidth}
    \centering
    \includegraphics[width=\textwidth]{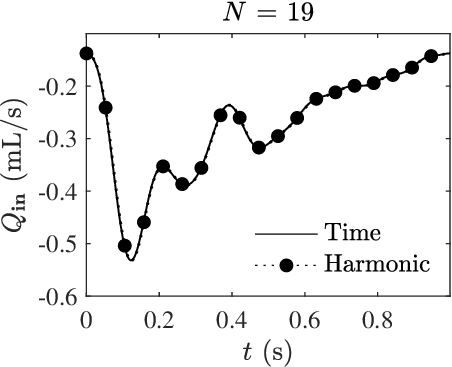}
    \caption[]{{\small}}
  \label{fig:cere_bc_n19}
  \end{subfigure}}
  \caption{Discretized harmonic balance boundary conditions at each time point (black circle) on the LVA and RVA inlets. The dotted line represents the reconstructed profile from the time points. For reference, the solid line is the original flow profile in time. (a) $N = 7$, (b) $N = 13$, (c) $N = 19$.}
  \label{fig:cere_bc}
\end{figure}

Four simulations are performed: one using the conventional time-stepping solver and three using the harmonic balance solver with different numbers of time points ($N = 7$, 13, and 19). 
Figure~\ref{fig:cere_bc} presents the harmonic balance boundary conditions. These flow profiles were prescribed identically at both the LVA and the RVA inlets. 
The geometry is discretized using 734,863 tetrahedral elements.
The pseudo time step size is $2 \times 10^{-2}$ resulting in a CFL number of $C_{\mathrm{CFL}} = 2$. 

\begin{figure} [h]
    \centering
    \includegraphics[width=0.5\linewidth]{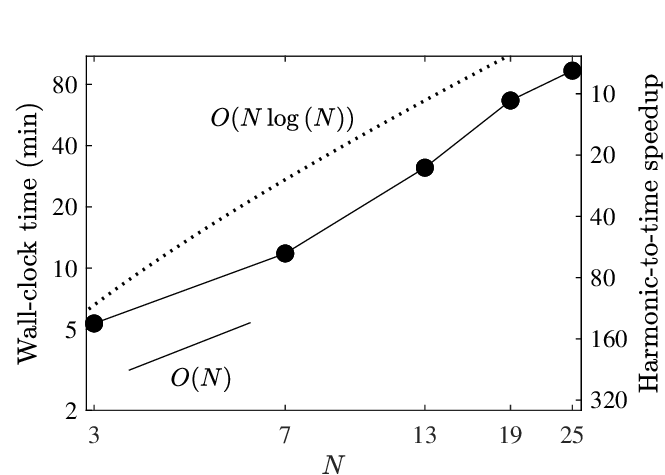}
    \caption{The cost performance of the harmonic balance solver for the cerebral flow case. The left axis shows the wall-clock time of the harmonic balance simulations. The right axis shows the speedup of the harmonic balance simulations compared to the conventional time simulation.}
    \label{fig:cere_cost}
\end{figure}

The simulation cost of the harmonic balance solver is shown in Figure~\ref{fig:cere_cost}. 
The conventional time simulation took 12 hours to run four cardiac cycles using 1,000 time steps per cycle.
The cases with $N = 3$ and 25 are simulated for demonstrating the cost scaling.
The cost scaling skews toward $O(N \log{(N)})$ as the number of time points, $N$, gets larger, which is as expected from the fast Fourier transforms implementation of the harmonic balance solver.

\begin{figure} [H]
    \centering
    \includegraphics[width=0.8\linewidth]{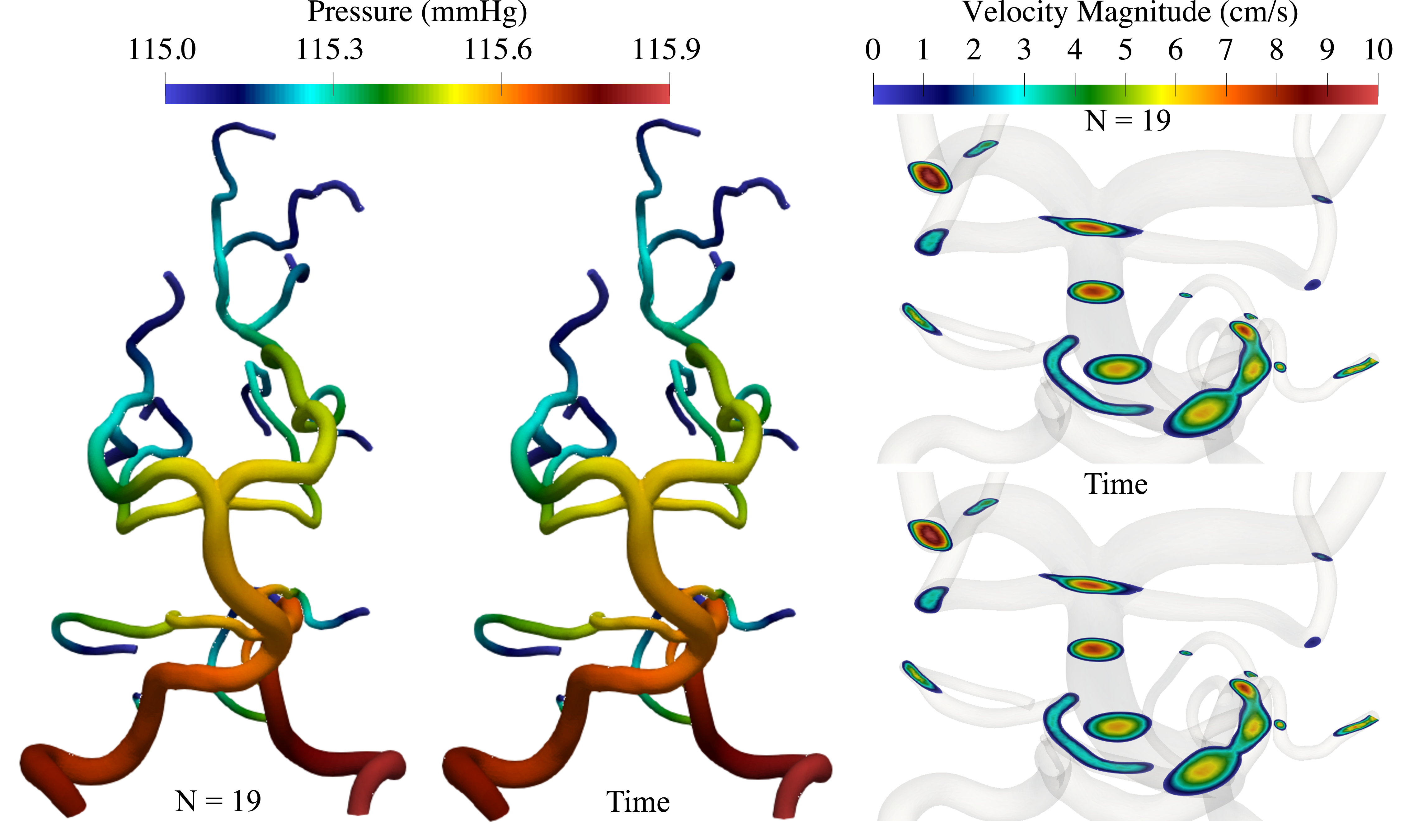}
    \caption{Harmonic balance results using $N = 19$ compared to the conventional time results for the pulmonary flow case. The pressure and velocity magnitude contours are taken at $t = 0.2$ seconds of the cardiac cycle.}
    \label{fig:cere_cont}
\end{figure}

The pressure and velocity contours for the harmonic balance solution using $N = 19$ and the conventional time solution at $t = 0.2 s$ in the cardiac cycle are shown in Figure~\ref{fig:cere_cont}. 
As shown in the figure, the results show no noticeable qualitative differences.

\begin{figure} [H]
    \centering
    \includegraphics[width=0.5\linewidth]{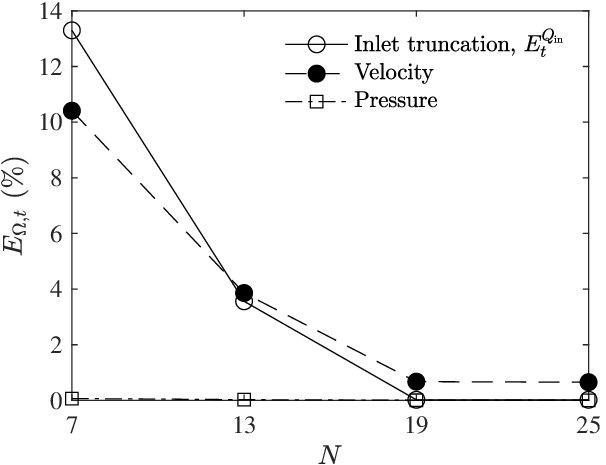}
    \caption{The volumetric and temporal integral relative root mean square error over the fluid domain and one cardiac cycle of the cerebral flow results.}
    \label{fig:cere_eint}
\end{figure}

Similar to the pulmonary flow case, we can calculate the integral relative errors.
Figure~\ref{fig:cere_eint} shows the spacial and time integral error of the solution along with the truncation error for the inlet boundary condition. 
In this case, due to the low Reynolds number, the inlet flow truncation error is a good indicator of the velocity field solution error.
Although we took conventional CFD as gold standard in our error calculations, the true error of harmonic balance solver (if the exact solution was available) may have been smaller given that the our ``gold standard'' solution itself contains some error. 
Similar to the earlier case, the pressure error is considerably smaller owing to the large fixed pressure boundary conditions at the outlets.

Figure~\ref{fig:cere_eint} includes an additional set of results at $N = 25$.
Since the inlet boundary profile is exactly represented with $N = 19$, any simulation run using more than 19 time points will not carry an inlet truncation error.
As the number of time points increases further past $N = 19$, the solution errors remain unchanged. 
This suggests that the reported error stems from the spatial discretization of both the conventional and harmonic balance formulations, rather than the frequency truncation of the solution. 
Determining which formulation more closely aligns with physiological results would require extensive clinical-to-simulation comparison studies, which are beyond the scope of this study.

\subsection{Left main coronary arteries flow} 

\begin{figure} [H]
    \centering
    \includegraphics[width=0.7\linewidth]{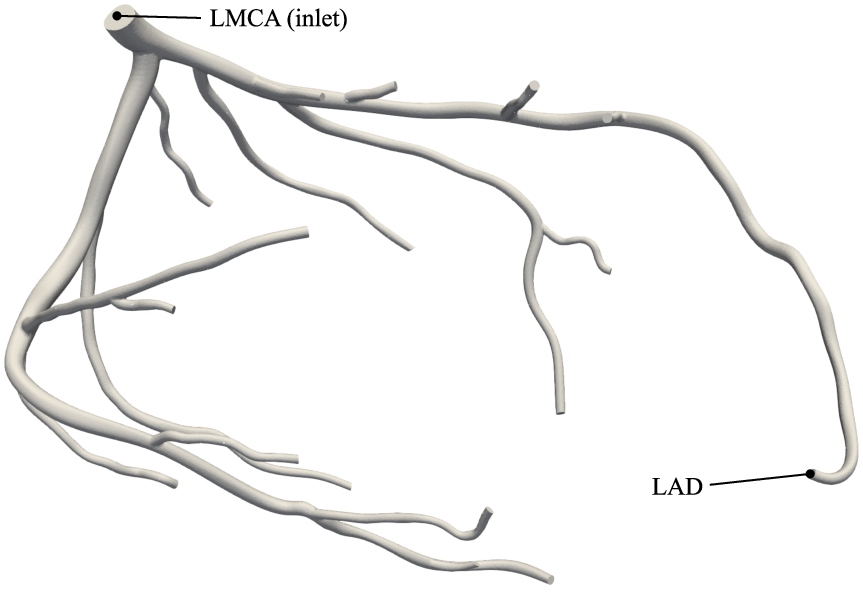}
    \caption{Clinically obtained left main coronary arteries (LMCA) geometry used in this study. The left anterior descending artery (LAD) is marked as the location for results in Figure~\ref{fig:lca_eq}. }
    \label{fig:lca_geo}
\end{figure}

Coronary artery disease is the leading cause of death in the United States~\cite{shao2020coronary,libby2005pathophysiology}.
It is predominantly caused by atherosclerosis of coronary arteries, leading to narrowing or blockage of the blood supply to the heart.
In the final case of this study, we will evaluate the proposed harmonic balance solver using a left main coronary arteries (LMCA) flow case. 
The LMCA originates from the ascending aorta and supplies blood to the left side of the heart, including the left ventricle. 
Computational fluid dynamics (CFD) simulations have been employed as a noninvasive method to gather flow rate and pressure information on patient-specific coronary geometries. 
This information can help determine the severity of coronary artery disease and assist in surgical decision-making.
The LMCA geometry used in this case was segmented from a clinically obtained image~\cite{wilson2013vascular}.
The aorta and right coronary arteries are removed from the original geometry so that only the LMCA remains, as shown in Figure~\ref{fig:lca_geo}.
This decision was made due to the flow instabilities and backflow introduced by the short aorta in the original image, which made simulation convergence impossible for both the conventional CFD solver and the proposed harmonic balance solver.
The inlet flow profile taken from literature~\cite{berne2001cardiovascular} was used as the boundary condition for prior conventional CFD studies~\cite{olgac2009patient}. 
This inlet boundary condition poses numerical challenges because the inflow profile contains distinctive kinks, which require a substantial number of time points to represent accurately using the harmonic balance solver.
However, we will demonstrate that the solver produces reasonable solutions using a manageable number of time points.
The Neumann outlets are set at a steady pressure of 70 mmHg.
The resulting peak Reynolds number of the domain is around 100.

\begin{figure}[H]
    \centering
{\captionsetup{position=bottom,justification=centering}\begin{subfigure}[b]{0.32\textwidth}
    \centering
    \includegraphics[width=\textwidth]{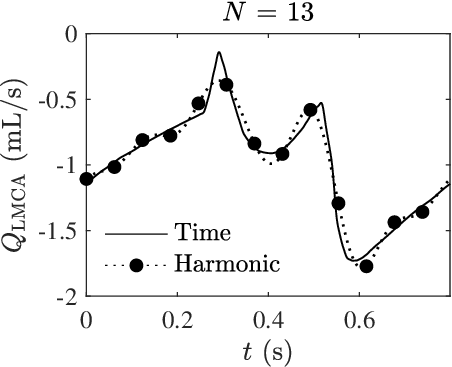}
    \caption[]{{\small}}
    \label{fig:lca_bc_n7}
  \end{subfigure}
\hfill  \begin{subfigure}[b]{0.32\textwidth}
    \centering
    \includegraphics[width=\textwidth]{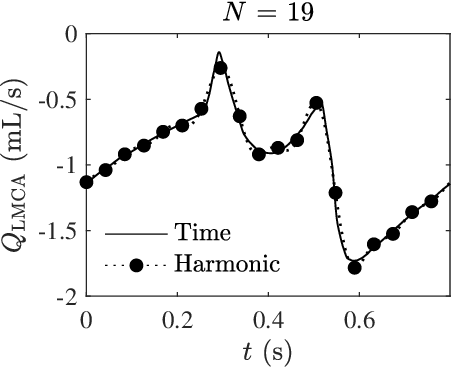}
    \caption[]{{\small}}
    \label{fig:lca_bc_n13}
  \end{subfigure}
\hfill
    \begin{subfigure}[b]{0.32\textwidth}
    \centering
    \includegraphics[width=\textwidth]{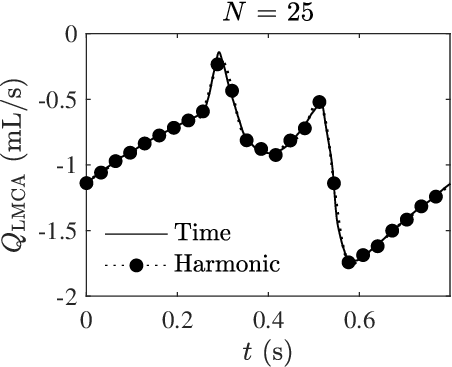}
    \caption[]{{\small}}
  \label{fig:lca_bc_n19}
  \end{subfigure}}
  \caption{Discretized harmonic balance boundary conditions at each time point (black circle) on LMCA inlet. The dotted line represents the reconstructed profile from the time points. For reference, the solid line is the original flow profile in time. (a) $N = 13$, (b) $N = 19$, (c) $N = 25$.}
  \label{fig:lca_bc}
\end{figure}

Six simulations were performed, one using the conventional time solver and five using the harmonic balance solver ($N = 1$, 7, 13, 19, and 25). 
We simulated a case using the steady-state assumption ($N=1$) to demonstrate that this assumption, that may be used to reduce computational costs in clinical settings, is insufficient for capturing the full dynamics of the flow.
The inlet boundary conditions for $N = 13$, 19, and 25 are shown in Figure~\ref{fig:lca_bc}. The period of the cardiac cycle is $T = 0.8$ seconds.
As observed, the sharp changes in the flow rate profile cannot be fully captured by 13 time points.
As the number of time points increases, the inflow profile is more accurately represented. 
The computational domain is discretized using 1,336,207 tetrahedral elements. 
The pseudo time step size is $5 \times 10^{-3}$, which results in a CFL number of $C_{\mathrm{CFL}} = 1.3$. 

\begin{figure} [H]
    \centering
    \includegraphics[width=0.5\linewidth]{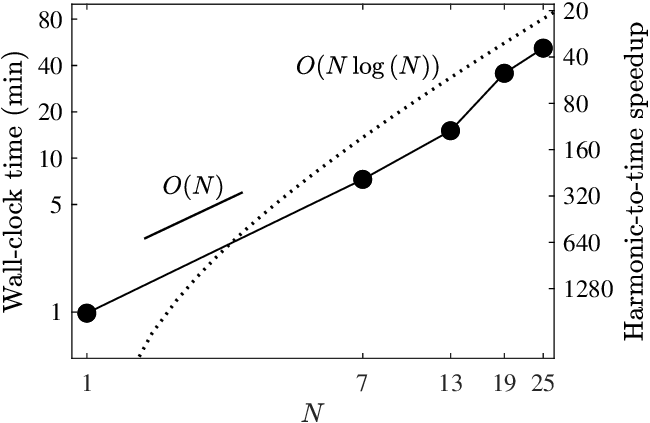}
    \caption{The cost performance of the harmonic balance solver for the coronary flow case. The left axis shows the wall-clock time of the harmonic balance simulations. The right axis shows the speedup of the harmonic balance simulations compared to the conventional time simulation.}
    \label{fig:lca_cost}
\end{figure}

The simulation wall clock time and speedup of the harmonic balance solver are summarized in Figure~\ref{fig:lca_cost}. 
The conventional time simulation took 30.3 hours to run for 4 cardiac cycles, using 800 time steps per cycle. 
The cost scaling is similar to that of the cerebral flow cases, going from $O(N)$ to $O(N\log{(N)})$ as $N$ increases. 
The harmonic balance solver is more than $30$ times faster than the conventional time solver, for the case using the largest number of time points ($N = 25$). 

\begin{figure} [H]
    \centering
    \includegraphics[width=0.7\linewidth]{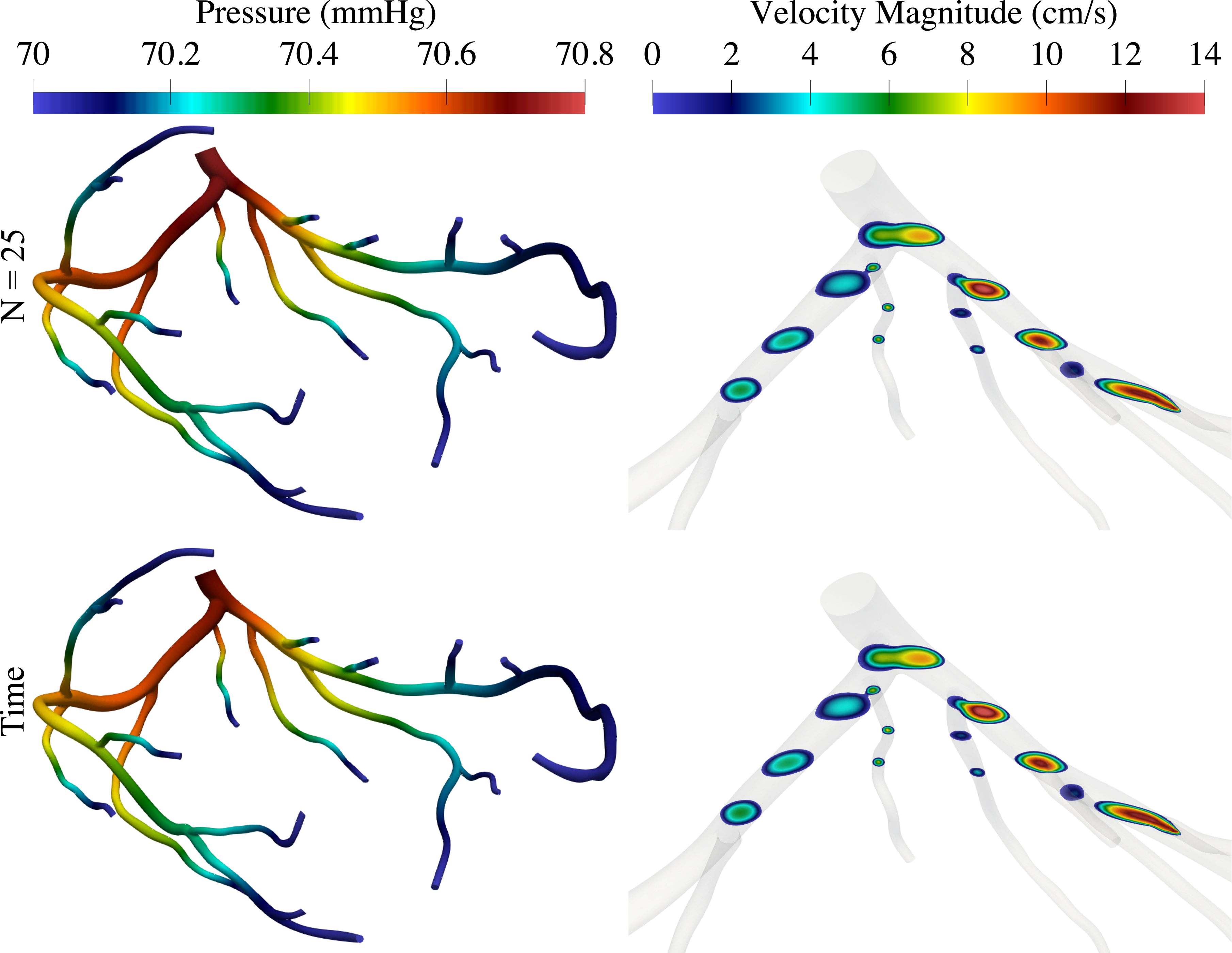}
    \caption{Harmonic balance results using $N = 25$ compared to the conventional time results for the coronary flow case. The pressure and velocity magnitude contours are taken at $t = 0.5$ seconds of the cardiac cycle.}
    \label{fig:lca_cont}
\end{figure}

The pressure and velocity magnitude contours for the harmonic balance solver solution using $N = 25$ and the conventional time solver solution are shown in Figure~\ref{fig:lca_cont}. 
The snapshots are taken at $t = 0.5 s$ seconds in the cycle. 

\begin{figure}[H]
{\captionsetup{position=bottom,justification=centering}\begin{subfigure}[b]{0.48\textwidth}
    \centering
    \includegraphics[width=\textwidth]{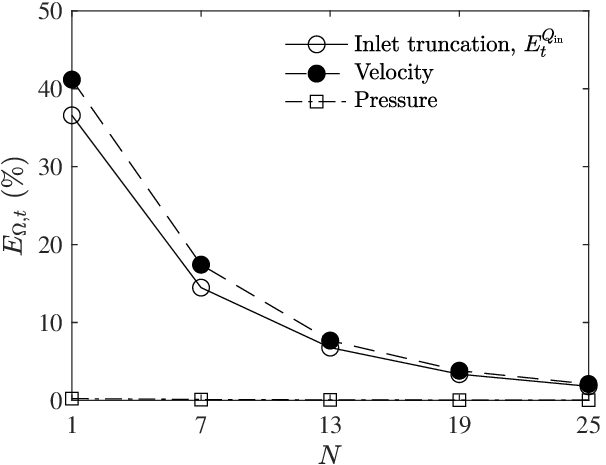}
    \caption[]{{\small}}
    \label{fig:lca_e_int}
  \end{subfigure}
  \hfill
  \begin{subfigure}[b]{0.48\textwidth}
    \centering
    \includegraphics[width=\textwidth]{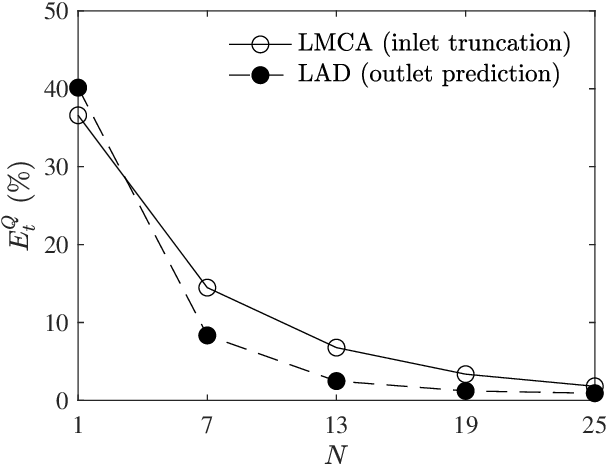}
    \caption[]{{\small}}
    \label{fig:lca_eq}
  \end{subfigure}}
  \caption{Left main coronary arteries flow relative root mean square error for (a) the velocity and pressure integrated over the fluid domain and one cardiac cycle and (b) the inlet LMCA flow rate (truncation) and the outlet LAD flow rate (prediction). }
  \label{fig:lca_eint}
\end{figure}

The integral relative error in Figure~\ref{fig:lca_e_int} shows that the harmonic balance velocity result is within $5\%$ of that of the conventional time solver results when $N \geq 19$. 
The steady state flow assumption ($N=1$) produces significant velocity error of more than $40\%$.
Similar to the cerebral flow case, the inlet boundary condition truncation error is a good estimate of the overall error in the solution domain due to the flow's low Reynolds number. 

\begin{figure}[H]
  {\captionsetup{position=bottom,justification=centering}\begin{subfigure}[b]{0.32\textwidth}
    \centering
    \includegraphics[width=\textwidth]{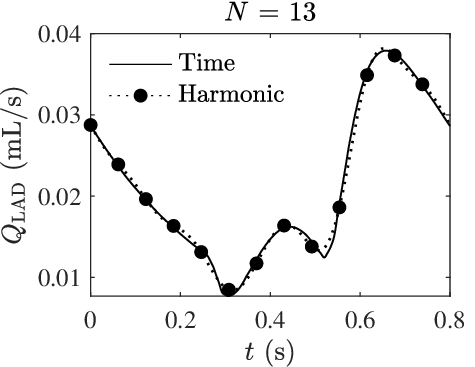}
    \caption[]{{\small}}
    \label{fig:lcda_n13}
  \end{subfigure}
  \hfill
  \begin{subfigure}[b]{0.32\textwidth}
    \centering
    \includegraphics[width=\textwidth]{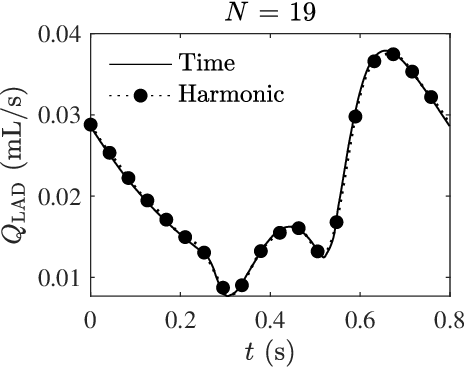}
    \caption[]{{\small}}
    \label{fig:lcda_n19}
  \end{subfigure}
\hfill
    \begin{subfigure}[b]{0.32\textwidth}
    \centering
    \includegraphics[width=\textwidth]{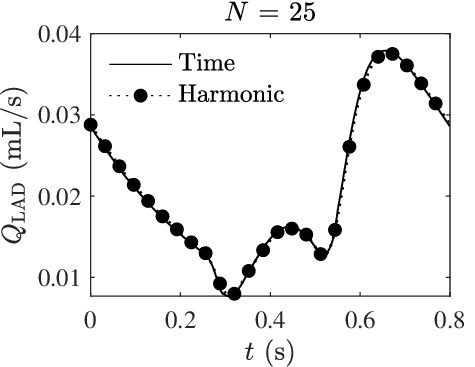}
    \caption[]{{\small}}
  \label{fig:lcda_n25}
  \end{subfigure}}
  \caption{left anterior descending artery (LAD) flow rate results obtained with the harmonic balance solver (solid circles) and the with the conventional time solver (solid line). The dotted line represents the reconstructed profile from the harmonic balance time points. (a) $N = 13$, (b) $N = 19$, (c) $N = 25$.}
  \label{fig:lad}
\end{figure}

For some parameters of interest, the error in the results will be smaller than the inlet truncation error due to the damping effect of viscosity, which smooths out sharp changes in velocity downstream of the flow. 
For instance, the flow rate at the left anterior descending artery (LAD), marked in Figure~\ref{fig:lca_geo}, can be accurately predicted using just 13 time points ($N=13$), as demonstrated in Figure~\ref{fig:lad}.
The LAD outlet flow rate relative error is under $3\%$ when $N \geq 13$.
Aside from the steady flow case, the flow rate errors at the LAD outlet are significantly lower than the inlet truncation errors at the left main coronary artery (LMCA).

\section{Discussion}
The above results demonstrated that the proposed harmonic balance CFD solver is a viable option for various cardiovascular simulations, providing a significant speedup over to the current time-stepping CFD solver. 
The requirements for conducting harmonic balance CFD simulations are close to those of conventional CFD simulations: a discretized computational domain (mesh) and boundary conditions.
The harmonic balance solver has two additional prerequisites regarding the flow characteristics: 1. the flow must exhibit periodic behavior with a known frequency, which is fundamentally satisfied by cardiovascular flows, and 2. there should be no geometry-triggered instabilities or turbulence in the fluid field, meaning that the Reynolds number has to be below the transitional flow Reynolds number, typically between 2,000 and 3,000. 
In cardiovascular flows, the Reynolds number can exceed this limit in the aorta, where the peak Reynolds number can reach 4,000. This is the reason why we excluded the aorta in the coronary flow case. 
Nonetheless, a significant portion of the human circulatory system operates within the acceptable limit, making the proposed harmonic balance solver beneficial for simuilating the flow in those areas. 
The simplicity and adaptability of the harmonic balance CFD solver set it apart from other cost-reducing methods, such as lumped parameter networks or data-driven models, which require extensive prior training and model development to yield satisfactory results.
Additionally, the harmonic balance solver is grounded in fundamental conservation laws, which other methods mentioned above do not guarantee.

In this study, the harmonic balance solver outperforms the conventional time solver by 10 to 100 times in simulation speed across all cases.
This remarkable speedup is achieved by implementing a fast Fourier transform scheme for better cost scaling as the number of time points increases. 
The improvement in cost scaling distinguishes this study from our earlier study, where the second-order cost scaling ($O(N^2)$) diminished the speedup advantage as $N$ becomes large~\cite{esmaily2024new}.
The harmonic balance solver retains the same parallel scalability benefits as the previous study~\cite{esmaily2024new}, which came for ``free'' without modifications to the solver's parallelization scheme or code structure.
\begin{figure}
    \centering
    \includegraphics[width=0.5\textwidth]{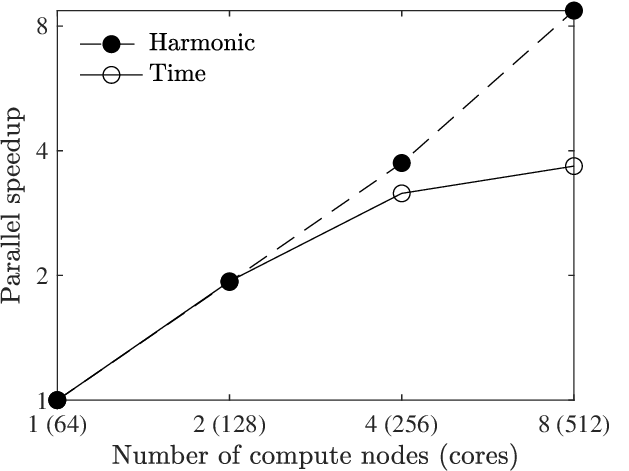}
    \caption{Strong scaling of the harmonic balance solver and the conventional time solver using the $N = 25$ case in Section 3.3.}
    \label{fig:npscale}
\end{figure}
As shown in Figure~\ref{fig:npscale}, the harmonic balance simulations can be run on more cores in parallel than the conventional time CFD solver while still maintaining strong scaling. 
This capability arises from a denser linear system at each mesh node, which results from solving for all the time points simultaneously. 
Consequently, this higher density helps maintain strong scaling, even when fewer mesh nodes are assigned to each parallel partition as the number of cores increases. 

\begin{figure}
    \centering
    \includegraphics[width=0.5\textwidth]{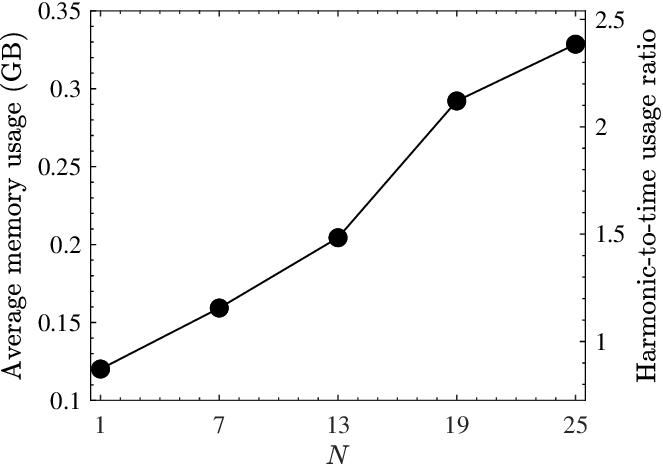}
    \caption{Memory usage for the harmonics balance solver and compared to the conventional solver for the LMCA cases in Section 3.3. The left axis shows the average memory usage over the simulations. The right axis shows the average memory usage ratio compared to the  conventional time-stepping solver.}
    \label{fig:mem}
\end{figure}

The harmonic balance solver computes all time points simultaneously, which can result in substantial memory usage—specifically, $O(N^2)$ times the memory required by a conventional time-stepping solver, if implemented with a brute-force full tangent matrix approach.
In this study, we significantly reduced memory usage by combining several techniques: matrix splitting, fast Fourier transformations, and an efficient sparse matrix multiplication scheme that condenses repetitive non-zero entries.
Figure~\ref{fig:mem} presents the average memory usage for all left main coronary artery cases discussed in Section 3.3. 
The average memory usage was calculated from dividing the total memory usage (in gigabyte-seconds) by the CPU time (in seconds), as reported directly by the grid engine scheduler.
For reference, the simulation ran with the conventional time-stepping solver has an average memory usage of 0.138 gigabytes. 
As shown in the figure, the harmonic balance solver simulating 25 time points ($N=25$) required only $2.4$ times the memory of the conventional solver. 
This marginal increase in memory footprint is a result of the fact that the bulk of memory in finite element simulations is dedicated to storing arrays other than unknowns, such as sparse matrices and connectivity information. 

The simulation results show that the truncation error from the boundary condition is a reliable indicator for the overall solution error when the flow Reynolds number is $O(100)$.
In these situations, he number of time points ($N$) chosen in the harmonic balance simulations can be determined by the acceptable level of solution error, which is reflected by the boundary condition truncation error.
For simulations with intermediate Reynolds numbers around or over $1,000$, additional time points are beneficial for accessing the higher-frequency flow information generated by the nonlinear convective term. 
In addition, the truncated frequency representation of the flow helps eliminate non-physical oscillations in the solution. 
Clinical flow measurements, often prescribed as boundary conditions, are noisy. 
The harmonic method acts as a built-in filter, effectively removing high-frequency signals and eliminating non-physical oscillations caused by the noise in the inflow signals.

It is important to emphasize that the study does not determine the actual error of the harmonic balance solver in relation to real-life physical measurement; all errors are reported relative to conventional time-stepping solvers.
Thus, the reported error can over-estimate the actual error since they are measured relative to the conventional method, the solution of which contains some error.
Our previous paper provided a comprehensive error analysis and established the bounds of the numerical errors associated with the frequency solver~\cite{esmaily2024new}. 

One potential improvement to the proposed solver is the implementation of resistance outlet boundary conditions. 
This boundary type is considered to be more physiologically accurate~\cite{troianowski2011three}. 
For resistance boundaries, the pressure at the boundary is proportional to the flow rate that passes through that boundary.
Consequently, the pressure error resulting from resistance boundaries is likely to be similar to the error of outlet flow rate.
Future studies are required to confirm this hypothesis.

Another potential application of the harmonic balance solver include respiratory flow simulations, which are also periodic in time, and the peak Reynolds number is typically less than 2,000. 
Broadly speaking, the harmonic balance CFD solver is an excellent tool for any other application with a time-periodic flow field and low to intermediate Reynolds number flow.
The use of harmonic balance CFD solver solutions to build lumped parameter networks or generate machine learning datasets can also be studied.
It has the potential to significantly speed up the process of constructing and training these models.
The proposed solver can also be coupled with fluid-structure interaction codes, in which the governing equations must be formulated in the harmonic balance form. 

\section*{Data availability}
The simulation files and results presented in this study are available on request from the corresponding author.

\section*{Acknowledgement}
For the patient-specific case in Section 3, the geometry used herein was provided in whole or in part with Federal funds from the National Library of Medicine under Grant No. R01LM013120, and the National Heart, Lung, and Blood Institute, National Institutes of Health, Department of Health and Human Services, under Contract No. HHSN268201100035C.

\section*{Author contributions statement}
Conceptualization, D.J. and M.E.; methodology, D.J. and M.E.; software, D.J. and M.E.; validation, D.J.; formal analysis, D.J.; investigation, D.J.; resources, M.E.; data curation, D.J.; writing---original draft preparation, D.J.; writing---review and editing, D.J. and M.E.; visualization, D.J.; supervision, M.E.; project administration, M.E.. All authors have read and agreed to the published version of the manuscript.

\bibliographystyle{unsrt}
\bibliography{ref}

\end{document}